\documentclass[amstex, 10pt]{article}
\usepackage[margin=1in]{geometry}
\usepackage{amsmath}
\usepackage{amsthm}
\usepackage{amssymb}
\usepackage{authblk}
\allowdisplaybreaks
\newtheorem{theorem}{Theorem}[section]
\newtheorem{lemma}[theorem]{Lemma}

\title{\bf Observability and Control Property for a Singular Heat Equation  with Variable Coefficients}
\author [a] {Xue Qin \thanks{ Email: \,\, qx3@mail.ustc.edu.cn. }}
\author [b] { Shumin Li \thanks{ Email:\,\, shuminli@ustc.edu.cn.}}
\affil [a] {School of Mathematical Sciences,
 University of Science and Technology of China, No.96, JinZhai Road, Hefei, Anhui, 230026, P.R.CHINA}
 \affil [b] {Wu Wen-Tsun Key Laboratory of Math, USTC, Chinese Academy of Sciences. School of Mathematical Sciences,
 University of Science and Technology of China, No.96, JinZhai Road, Hefei, Anhui, 230026, P.R.CHINA }

\date{}

\renewcommand{\baselinestretch}{1.0}

\begin{document}
\maketitle

 \begin{abstract}
 The goal of this paper is to analyze control properties of the parabolic equation with variable coefficients in the principal part and with a singular inverse-square potential:\,$\partial_tu(x,t)-{\rm div}(p(x)\nabla u(x,t))-({\mu}/{|x|^2})u(x,t)=f(x,t).$  Here $\mu$ is a real constant . It was proved in the paper of Goldstein and Zhang (2003) that the equation is well-posedness when $0\leq{\mu\leq p_1(n-2)^2/4}$, and in this paper, we mainly consider  the case $0\leq\mu<({ p_1^2}/{ p_2})(n-2)^2/4$ , where $ p_1,p_2$ are two positive constants which satisfy:\, $ 0< p_1\leq p(x)\leq p_2 , \forall x\in\overline{\Omega}$. We extend the specific  Carleman estimates in the paper of Ervedoza  (2008) and Vancostenoble (2011) to the equation we consider and apply it to deduce an observability inequality for the system. By this inequality and the classical HUM method, we obtain that we can control the equation from any non-empty open subset as for the heat equation. Moreover, we will study the case $\mu>p_2(n-2)^2/4$. We consider  a sequence of regularized potentials $\mu/(|x|^2+\epsilon^2),$ and prove that we cannot stabilize the corresponding systems uniformly with respect to $\epsilon>0,$ due to the presence of explosive modes which concentrate around the singularity.

  {\bf MSC} 35B37;\,93B05;\,35K05;\,49J20.

  {\bf Keywords}  Singular heat equation; Null controllability; Carleman estimate; Observability.
 \end{abstract}

\section{Introduction and main results}
\setcounter{equation}{0}

The purpose of this work is to establish null control for linear heat equations with variable
coefficients in the principal part and with singular potentials. We are interested in the so-called inverse square potential of  the form $-\mu/|x|^2$ which arises, for example, in combustion theory \cite{BJ,BV,D,GV} and quantum mechanics \cite{BG2,d,RM}.

 Let $n\geq3$ and consider a  bounded domain $\Omega\subset \mathbb{R}^n$ which satisfies $0\in\Omega$ and the boundary $\partial \Omega\in C^2$.
$\omega\subset\Omega$ is a non-empty open set which can be fixed arbitrarily.

 We are interested in the control and stabilization properties of the following system
\begin{equation}\label{1.1}
 \begin{cases}
 \partial_tu(x,t)+Lu(x,t)=f(x,t), \quad{(x,t)\in\Omega_T,}\\
 u(x,t)=0,\qquad{ (x,t)\in\Sigma_T,}\\
 u(x,0)=u_0(x),\qquad {x\in\Omega.}
 \end{cases}
 \end{equation}
Here, 
\begin{gather}
Lu(x,t)=-{\rm div}(p(x)\nabla u(x,t))-\frac{\mu}{|x|^2}u(x,t),
\end{gather}
 $x=(x_1,\ldots,x_n)\in \Omega\subseteq \mathbb{R}^n,$ $\partial_t=\frac{\partial}{\partial_t}$, $\partial_j=\frac{\partial}{\partial x_j}, j=1,\ldots,n$, $\nabla=(\partial_1,\ldots,\partial_n)$. $\mu$ is a real constant, $u_0\in L^2(\Omega)$, and
  \begin{align}\label{Ome}
  \small\Omega_T=\Omega\times(0,T), \qquad \small\omega_T=\omega\times(0,T),\qquad \small\Sigma_T=\partial\Omega\times(0,T).
  \end{align}
 Assume that
 \begin{align}\label{a}
 p(x)\in C^3(\overline{{\Omega}}),\,\,
 0<p_1\leq p(x)\leq p_2,\,\,|\nabla p|\leq p_3, \qquad x\in\overline{{\Omega}},
 \end{align}
 for some constants $0<p_1\leq p_2$ and $0\leq p_3$.
 $f\in L^2((0,T);H^{-1}(\Omega))$ is the control that we assume to be null in $\Omega\backslash\bar{{\omega}}$, that is
\begin{gather}\label{*}
 \forall \theta\in \mathcal{D}(\small\Omega\backslash\bar{{\omega}}),\quad \theta f=0 \quad in \quad L^2((0,T);H^{-1}(\Omega)).
\end{gather}

First of all, let us briefly mention that the Cauchy problem with such singular potential $-\mu/|x|^2$ is not straightforward. Indeed, Goidstein and Zhang \cite{GZ1} proved that the existence and nonexistence of positive solutions to (1.1) is largely determined by the size of the infimum of the spectrum of the symmetric operator with more general singular potentials $V$: $L=-\sum\limits_{i,j=1}^n\partial_j(a_{ij}(x)\partial_iu(x,t))-V(x)u(x,t),$  and $0\leq V\in L_{loc}^1(\Omega).$ Here, $\{a_{ij}\}_{1\leq i\leq j\leq n}$ satisfy  $a_{ij}=a_{ji}$ and  the uniform ellipticity. Actually, this problem is also strongly related to the Hardy inequality:
 \begin{equation}\label{1.3}
 \mu^*(n)\int_\Omega \frac{u^2}{|x|^2}\,{\rm d}x\leq\int_\Omega |\nabla u|^2\,{\rm d}x,\qquad  \forall u\in H_0^1(\Omega),
  \end{equation}
  where $\mu^*\triangleq\mu^*(n)=(n-2)^2/4$ is the optimal constant. Using the principle of Duhamel and the direct applications of  Theorem 1.1 (iii) in \cite{GZ1}, (\ref{a}), and (\ref{1.3}), respectively, we  can obtain:
if $0\leq\mu\leq p_1\mu^*$, the system (1.1) has a unique nonnegative solution  $u(\cdot, t)\in W_{\rm loc}^{1,2}(\Omega) $ for any $0\leq u_0\in L^2(\Omega),$ and $f(\cdot,t)\in L^2(\Omega)$ for all $t\in (0,T)$.

The first work about the existence and nonexistence of positive solutions to (1.1) with $p(x)=1$ was discovered by Baras and Goldstein \cite{BG1} in 1984. It was proved that if $0\leq\mu\leq\mu^*$ and if the initial data $u_0$ is positive,  it has a global weak solution , and  if $\mu>\mu^*,$ it has no solution if $u_0>0$ and $f\geq0$, even locally in time. The proof in \cite{BG1} relies on Hardy's inequality and scaling properties of the heat equation in $\mathbb{R}^n$ and, \cite{BG1} does not give relevant results about parabolic equations with variable coefficients in the principal part. In 1999, the case  $V\in L^1_{loc}(\Omega)$ was settled by Cabr\'{e} Xavier and Martel Yvan (\cite{CM}). The proof in \cite{CM} involves nice and delicate computations using the special structure of $\Delta$ and the smoothness of the boundary. In 2003, using an extension of the method in \cite{CM}, Goidstein and Zhang \cite{GZ1} extended the existence and blow up result  by Baras and Goldstein \cite{BG1}  to parabolic equations with variable
leading coefficients under almost optimal conditions on the singular potentials.
  In fact, the proof of \cite{GZ1} is greatly influenced by the ideals of \cite{CM} and the Harnack chin arguments of \cite{GZ2}.  Moreover, the condition of \cite{GZ1} on the singular potential allows new nonradial singularities such as $\mu(x)/|x|^2$ with $\mu(x)$ being an unbounded function and potentials with nonpointwise singularities such as $c/{d(x,\partial \Omega)}^2.$

{\it Null-Controllability problem:

   Given any $u_0\in L^2(\Omega),$ find a function $f\in L^2(\omega\times(0,T))$ satisfying $(\ref{*})$  such that the solution of $(1.1)$ satisfies
  \begin{gather}\label{null}
    u(x,T)=0,\qquad x\in \Omega.
  \end{gather}}

Notice that for $\forall 0\leq\mu\leq p_1\mu^*,$  we only need to consider the case $0\notin \omega$. In fact, for the case $0\in\omega$, the stabilization property of (1.1) is always holds for all $\mu$ (see Section 5).

 A special controllability issue was already discussed under the assumption $0\leq\mu\leq\mu^*$ when $p(x)=1$  in the recent work \cite{VZ1}, in the special case where $\omega$ contains an annulus centered in the singularity. The authors of \cite{VZ1} need this assumption since their proof strongly uses a decomposition in spherical harmonics which allows  to reduce the problem to the study of 1-d singular equations. Later on, \cite{E} proved that we can actually remove this assumption. It extended the 1-d Carleman estimates to the N-d case and deduced a null controllability result for a control supported in any non-empty subset $\omega$, still provided that $0\leq\mu \leq\mu^*$. Actually, the proof of current paper is greatly influenced by the ideals of \cite{E} and the inverse source problems of \cite{J} with shaper Carleman estimates.

In this article we proposed a study of a parabolic equation with variable coefficients in the principal part, from a control point of view, in the case $\mu < ({p_1^2}/{p_2})\mu^*(n)$, and the case $\mu>p_2\mu^*$. When $0\leq\mu< ({p_1^2}/{p_2})\mu^*(n)$, we have addressed the null-controllability problem for a distributed control in an arbitrary open subset of $\Omega$ (see Section 2). When $\mu>p_2\mu^*$, we have shown that we cannot uniformly stabilize regularized approximations of system (1.1) with a control supported in $\omega$ when $0\notin\omega$ (see Section 4).

 Now state the main result:
 \begin{theorem}
 \label{th1.1}
 Let $\mu$ be a real number such that $0\leq\mu<\frac{p_1^2}{p_2}\mu^*.$

 Given any non-empty open set $0\notin\omega\subset\Omega$, for any $T>0$ and $u_0\in L^2(\Omega)$, there exists a control $f\in L^2(\omega\times(0,T))$ such that the solution of (1.1) satisfies (\ref{null}). Besides, there exists a constant $C=C(T)$ such that
 \begin{gather}\label{control}
    \Vert{f}\Vert_{L^2(\omega\times(0,T))}\leq C(T) \Vert{u_0}\Vert_{L^2(\Omega)}.
  \end{gather}
  \end{theorem}
 Following the by now classical HUM method (\cite{HUM}), the controllability property is equivalent to an observability inequality for the adjoint system
 \begin{equation}\label{1.6}
  \begin{cases}
  \partial_ty(x,t)+{\rm div}(p(x)\nabla y(x,t))+\frac{\mu}{|x|^2}y(x,t)=0, \qquad {(x,t)\in\Omega_T,}\\
 y(x,t)=0, \qquad{ (x,t)\in\Sigma_T,}\\
 y(x,T)=y_T(x),\qquad {x\in\Omega.}
 \end{cases}
 \end{equation}
 More precisely, when $0\leq\mu<\frac{p_1^2}{p_2}\mu^*$, we need to prove that there exists a constant C such that for all $y_T\in L^2(\Omega)$, the solution of (\ref{1.6}) satisfies
  \begin{gather}\label{1.7}
  \int_\Omega\vert{y(x,0)}\vert^2\,{\rm d}x\leq C\int_{\omega_T}\vert{y(x,t)}\vert^2\,{\rm d}x\,{\rm d}t.
 \end{gather}

 In order to prove (\ref{1.7}), we will use a particular Carleman estimate, which is by now a classical technique in control theory, see for instance \cite{BD,FG2,FZ1,FZ2}. Indeed, the Carleman estimate we will derive later implies that for any solution $y$ of (\ref{1.6})
  \begin{gather}\label{1.8}
  \int_{\Omega\times(\frac{T}{4},\frac{3T}{4})}\vert{y(x,t)}\vert^2\,{\rm d}x\,{\rm d}t\leq C\int_{\omega_T}\vert{y(x,t)}\vert^2\,{\rm d}x\,{\rm d}t,
 \end{gather}
 which directly implies inequality (\ref{1.7}) since $t\mapsto\Vert{y(t,\cdot)}\Vert^2_{L^2(\Omega)}$ is increasing by the Hardy inequality (\ref{1.3}).
 In fact, when $\mu\leq p_1\,\mu^*$, multiplying (\ref{1.6}) by $y$, one obtains
 \begin{align}
 &\nonumber\frac{{\rm d}}{{\rm d}t}\int_\Omega |y(t,x)|^2\,{\rm d}x=-2\int_\Omega y\,{\rm div}(p(x)\nabla y)\,{\rm d}x-2\int_\Omega\frac{\mu}{|x|^2}|y(t,x)|^2\,{\rm d}x\\ \nonumber
 &=0+2\int_\Omega p(x)|\nabla y(t,x)|^2\,{\rm d}x
 -2\int_\Omega\frac{\mu}{|x|^2}|y(t,x)|^2\,{\rm d}x\\ 
 &\geq 2\int_\Omega \left(p_1|\nabla y(t,x)|^2-\frac{\mu}{|x|^2}|y(t,x)|^2\right)\,{\rm d}x\geq0.
  \end{align}

The Carleman estimate derived here is inspired by the works \cite{CM,MV} on 1-d degenerate heat equations, the recent paper \cite{VZ1} which is inspired  from the methods and results in \cite{CM,MV} to obtain radial estimates. As in \cite{CM,FI,MV,VZ1}, the major difficulty is to choose a special weight function appearing in the Carleman estimate. We adapt the same weight function with \cite{E,J,VZ1}.

This paper is organized as follows. In Section 2, we show the sharp Carleman estimate. In Section 3.1, we give the proof of the Carleman estimate. In Section 3.2, we give some proofs of technical lemmas. In Section 4, we will prove that we cannot  stabilize system (1.1) when $\mu> p_2\mu^*$ and $0\notin \omega$. In Section 5, we give some conclusions of this paper and point out that we can  stabilize system (1.1) when $0\in \omega$, for all $\mu$.

 \section{Null Controllability when $0\leq\mu< \frac{p_1^2}{p_2}\mu^*$}
 \setcounter{equation}{0}
In this section, we assume that $0\notin \omega.$  Moreover, for $\forall 0\leq\mu<\frac{p_1^2}{p_2}\mu^*$, we can take $r$ is  a fixed small positive constant which satisfies:
  \begin{gather}\label{r}
0<r<1,\quad 2p_1^2-2p_2\frac{\mu}{\mu^*}>3p_1p_3r,\\  \label{Br}
\overline B(0,r)\subseteq\Omega,\qquad \overline B(0,r)\cap\overline{\omega}=\varnothing.
\end{gather}
  Denote
 $B_r\triangleq B(0,r),$ $\overline B_r\triangleq\overline {B(0,r)}.$

 \subsection{Carleman Estimate}
 As said in the introduction, the main tool we use to address the observability inequality (\ref{1.8}) is a Carleman estimate. However, since it is based on tedious computations, we postpone the proofs of several technical lemmas to Section 3.

 The major problem when designing a Carleman estimate is the choice of a smooth weight function $\sigma$, which is in general assumed to be positive, and to blow up as t goes to zero and as t goes to T. That is:
 \begin{equation}\label{condition}
  \begin{cases}
  \sigma(x,t)>0, \qquad {(x,t)\in\Omega_T,}\\
 \lim\limits_{t \to 0^+}\sigma(x,t)=\lim\limits_{t \to T^-}\sigma(x,t)=+\infty,\qquad {x\in\Omega.}
 \end{cases}
 \end{equation}
 We propose the weight
 \begin{gather}\label{sigma}
 \sigma(x,t)=s\theta(t)\bigg({\rm e}^{2\lambda \sup\psi}-\frac{1}{2}|x|^2-{\rm e}^{\lambda\psi(x)}\bigg),
 \end{gather}
 where s and $\lambda$ are positive parameters aimed at being large,
  and $\psi$ is a smooth function satisfying
 \begin{equation}\label{psi}
  \begin{cases}
  \psi(x)=\ln\left(\frac{|x|}{r}\right), \qquad x\in B_r,\\
 \psi(x)=0, \,\,\,\,\,\,\,\,\,\,\,\,\,\qquad x\in\partial\Omega,\\
 \psi(x)>0,\,\,\,\,\,\,\,\,\,\,\,\,\,\qquad x\in\Omega\backslash\overline{B}_r,
 \end{cases}
 \end{equation}
 and there exists an open set $\omega_0$ and $\delta>0$ such that
 \begin{gather}\label{bounded}
 \overline{\omega}_0\subset\subset\omega,\quad |\nabla\psi(x)|\geq\delta,\quad x\in\overline{\Omega}\backslash\omega_0.
 \end{gather}
 The existence of such function $\psi$ is not straightforward but can be easily deduced from the construction given in \cite{FI}.

 Finally, $\theta$ is taken as \cite{J,VZ1}, in the following form:
 \begin{gather}\label{theta}
 \theta(t)=\bigg(\frac{1}{t(T-t)}\bigg)^{k_0},
 \end{gather}
 with $k_0>0$ to be chosen. Indeed, this is the key point that lead a shaper Carleman estimate.

 Note that the weight function $\sigma$ defined by (\ref{sigma}) indeed satisfies (\ref{condition}) and is smooth when $\lambda$ is large enough.

 To simplify notations, let us denote by $\phi$ the function
 \begin{gather}\label{phi}
 \phi(x)={\rm e}^{\lambda\psi(x)}
 \end{gather}
 and denote
  \begin{gather*}
  \Theta\triangleq\Omega\backslash\left(\overline{B}_r\cup\overline{\omega}_0\right),\quad \widetilde{\Theta}\triangleq\Omega\backslash\overline{B}_r,\\
 \Theta_T\triangleq\Theta\times(0,T),\quad \widetilde{\Theta}_T\triangleq\widetilde{\Theta}\times(0,T),\quad \omega_{0,T}\triangleq\omega_0\times(0,T).
 \end{gather*}

 We are now in position to state the Carleman estimate.
 \begin{theorem}
 \label{th2.1}
 Assume $0\leq\mu<\frac{p_1^2}{p_2}\mu^*$. Let $\gamma$ be given such that $0<\gamma<2$. Let $\sigma,\theta,\psi$ and $\phi$ be defined as above with
  $$k_0=1+2/\gamma.$$
 There exist some positive constants $C=C(\Omega,\omega,T,\gamma,p_1,p_2,p_3,r,\mu)>0,$ independent of $\lambda_0\geq1$ such that, for all $\lambda\geq\lambda_0$, there exists $s_0(\lambda)$ such that for all $s\geq s_0(\lambda)$, the solutions of (\ref{pm}) satisfy
\begin{align}\label{CE}
\nonumber&\,\,\,\,\,\,s\lambda^2\int_{\widetilde{\Theta}_T}\theta\phi {\rm e}^{-2\sigma}|\nabla w|^2\,{\rm d}x\,{\rm d}t
+\frac{s}{\lambda^2}\int_{\Omega_T}\theta {\rm e}^{-2\sigma}|\nabla w|^2\,{\rm d}x\,{\rm d}t\\ \nonumber
&+s\int_{\Omega_T}\theta {\rm e}^{-2\sigma}\frac{w^2}{|x|^\gamma}\,{\rm d}x\,{\rm d}t
+s^3\int_{\Omega_T}\theta^3{\rm e}^{-2\sigma}|x|^2w^2\,{\rm d}x\,{\rm d}t
+s^3\lambda^4\int_{\widetilde{\Theta}_T}\theta^3\phi^3{\rm e}^{-2\sigma}w^2\,{\rm d}x\,{\rm d}t\\ \nonumber
&+\frac{1}{s}{\rm e}^{-4\lambda\sup\psi}\int_{{B_r\times (0,T)}}\frac{1}{\theta}{\rm e}^{-2\sigma}w_t^2\,{\rm d}x\,{\rm d}t+\frac{1}{s}{\rm e}^{-4\lambda\sup\psi}\int_{\widetilde{\Theta}_T}\frac{1}{\theta\phi}{\rm e}^{-2\sigma}w_t^2\,{\rm d}x\,{\rm d}t\\ 
&\leq Cs^3\lambda^4{\rm e}^{2\lambda\sup\psi}\int_{\omega_T}\theta^3\phi^3{\rm e}^{-2\sigma}w^2\,{\rm d}x\,{\rm d}t
+ C\int_{\Omega_T}{\rm e}^{-2\sigma}g^2\,{\rm d}x\,{\rm d}t.
\end{align}
  Here $w$ satisfies:
  \begin{equation}\label{pm}
  \begin{cases}
  \pm\partial_tw(x,t)+{\rm div}(p(x)\nabla w(x,t))+\frac{\mu}{|x|^2}w(x,t)=g(x,t), \qquad {(x,t)\in\Omega_T,}\\
 w(x,t)=0, \qquad{ (x,t)\in\Sigma_T,}\\
 w(x,T)=w_T(x),\qquad {x\in\Omega,}
 \end{cases}
 \end{equation}
where, $w_T\in L^2(\Omega),\,g\in L^2(\Omega_T)$.
 \end{theorem}

The proof of Theorem 2.1 is given in Section 3.

\subsection{From the Caleman estimate to the Observability inequality}
In this Subsection, we explain why the Carleman estimate (\ref{CE}) implies the observability inequality (\ref{1.8}).

We apply Theorem to (\ref{1.6}),
so that we can  fix $\lambda>\lambda_0$ and $s>s_0$ such that (\ref{CE}) holds with $w=y$ and $g(x,t) \equiv 0$. These parameters now enter in the constant $C$:
\begin{align}\label{2.24}
\int_{\Omega_T}\theta {\rm e}^{-2\sigma}\frac{y^2}{|x|^\gamma}\,{\rm d}x\,{\rm d}t \leq C \int_{\omega_T}\theta^3\phi^3{\rm e}^{-2\sigma}y^2\,{\rm d}x\,{\rm d}t.
\end{align}
Noting the existence of constants $C$ such that
\begin{gather*}
\theta {\rm e}^{-2\sigma}\frac{1}{|x|^\gamma}\geq C,\qquad (x,t)\in \Omega \times\left[\frac{T}{4},\frac{3T}{4}\right],\\
\theta^3\phi^3{\rm e}^{-2\sigma}\leq C,\qquad (x,t)\in\omega\times(0,T).
\end{gather*}
We see that, (\ref{2.24}) implies (\ref{1.8}).

\section{Proof of Carleman Estimates: Proof of Theorem 2.1}
 \setcounter{equation}{0}
\subsection{Proof of Theorem 2.1}
Now we prove Theorem \ref{th2.1}.

 {\it 3.1.a. Notations and Preliminary Computations.}

 We present the main ideals and steps of the proof of Theorem 2.1, using several technical lemmas, that are proved later in Subsection 3.2.

 As usual in the proofs of Carleman inequalities, we define
 \begin{gather}\label{2.8}
z(x,t)=\exp(-\sigma(x,t))w(x,t),
\end{gather}
which obviously satisfies
\begin{gather}\label{boundary1}
z(T)=z(0)=0, \qquad in \quad H_0^1(\Omega)
\end{gather}
due to the assumptions (\ref{condition}) on $\sigma$.

Then, plugging $w=z\exp(\sigma(x,t))$ in the equation (\ref{pm}), we obtain that $z$ satisfies
\begin{align}\label{SP}
\nonumber\pm\partial_tz&+{\rm div}(p(x)\nabla z)+p\,|\nabla \sigma|^2z
+2\,p\,\nabla\sigma\cdot\nabla z\\ 
&+(\triangle\sigma)pz
+(\nabla p\cdot\nabla\sigma )z\pm (\partial_t\sigma) z+\frac{\mu}{|x|^2}z=g{\rm e}^{-\sigma}, \quad (x,t)\in\Omega_T,
\end{align}
with the boundary condition
\begin{gather}\label{boundary2}
z=0, \qquad (x,t)\in \Sigma_T.
\end{gather}

Following \cite{E}, since $r$ is a fixed number, we define a smooth positive radial function $\alpha(x)=\alpha(|x|)$ such that
\begin{equation}\label{alpha}
\begin{aligned}
\alpha(x)=0,\quad |x|\leq\frac{r}{2};\quad \alpha(x)=\frac{1}{n},\quad |x|\geq\frac{3r}{4},\\ 
0\leq\alpha(x)\leq\frac{1}{n},\quad \frac{r}{2}\leq|x|\leq\frac{3r}{4}.
\end{aligned}
\end{equation}

Setting the operators $S$ and $A$ are defined as below:
\begin{equation}\label{2.13}
\begin{aligned}
&Sz={\rm div}(p(x)\nabla z)+p\,|\nabla \sigma|^2z\pm(\partial_t\sigma)z+\frac{\mu}{|x|^2}z,\\ 
&Az=\pm\partial_tz+2\,p\,\nabla\sigma\cdot\nabla z+(1+\alpha)p(\triangle
\sigma)z+(\nabla p\cdot \nabla \sigma)z .
\end{aligned}
\end{equation}
One easily deduces from (\ref{SP}) that
\begin{gather*}
Sz+Az=\alpha\, p(\triangle\sigma)z+{\rm e}^{-\sigma}g,\\  \nonumber
 \Vert{Sz}\Vert^2+\Vert{Az}\Vert^2+2\left<Sz,Az\right>
\leq 2\left\Vert{\alpha \,p(\triangle\sigma) z}\right\Vert^2+2\Vert{{\rm e}^{-\sigma}g}\Vert^2,
\end{gather*}
where $\Vert{\cdot}\Vert$ denotes the $L^2(\Omega\times(0,T))$ norm and $\left<\cdot,\cdot\right>$  denotes the corresponding scalar product. Especially, the inequality
\begin{gather}\label{2.14}
\frac{1}{2}\Vert{Az}\Vert^2+\left<Sz,Az\right>-\left\Vert{\alpha \,p(\triangle\sigma)z}\right\Vert^2\leq
\Vert{{\rm e}^{-\sigma}g}\Vert^2,
\end{gather}
holds. Denote
\begin{gather}\label{I.0}
I=\left<Sz,Az\right>-\left\Vert{\alpha \,p(\triangle\sigma)z}\right\Vert^2.
\end{gather}

\begin{lemma}
 The following equality holds:
 \begin{align}\label{I}
 \nonumber I&=\int_{\Sigma_T}p^2(\nabla\sigma\cdot\vec\nu)|\partial_\nu z|^2{\rm d}S\,{\rm d}t
 +\int_{\Omega_T}\left[p\,(\nabla p\cdot\nabla\sigma)-\alpha \,p^2\triangle\sigma\right]|\nabla z|^2\,{\rm d}x\,{\rm d}t\\ \nonumber
 &-2\int_{\Omega_T}p(\nabla\sigma\cdot\nabla z)(\nabla p\cdot\nabla z)\,{\rm d}x\,{\rm d}t-2\int_{\Omega_T}p^2\sum_{k,j=1}^n(\partial_k\partial_j\sigma)(\partial_kz)\partial_jz\,{\rm d}x\,{\rm d}t\\ \nonumber
 &+\int_{\Omega_T}\left[\frac{1}{4}\bigg(3\nabla (p^2)\cdot\nabla\alpha+(1+\alpha)\triangle(p^2)\bigg)+\frac{1}{2}p^2\triangle\alpha\right]\triangle\sigma\,z^2\,{\rm d}x\,{\rm d}t\\ \nonumber
 &+\int_{\Omega_T}\left\{\left(\frac{3}{4}(1+\alpha)\nabla(p^2)+p^2\nabla\alpha\right)\cdot\nabla(\triangle\sigma)\right\}z^2\,{\rm d}x\,{\rm d}t\\ \nonumber
 &+\frac{1}{2}\int_{\Omega_T}\left[(1+\alpha)p^2\triangle^2\sigma+ \nabla p\cdot\nabla(\nabla p\cdot\nabla\sigma)+p\,\triangle(\nabla p\cdot\nabla\sigma)\right]\, z^2\,{\rm d}x\,{\rm d}t\\ \nonumber
 &+2\int_{\Omega_T}p\,(\nabla\sigma\cdot x)\frac{\mu}{|x|^4}\,\, z^2\,{\rm d}x\,{\rm d}t
 +\int_{\Omega_T}\alpha\, p(\triangle\sigma)\frac{\mu}{|x|^2}\,\, z^2\,{\rm d}x\,{\rm d}t-\frac{1}{2}\int_{\Omega_T}\left(\partial_t^2\sigma\right) z^2\,{\rm d}x\,{\rm d}t\\ \nonumber
 &\mp2\int_{\Omega_T} p\,(\nabla\sigma\cdot\nabla\partial_t\sigma)\, z^2\,{\rm d}x\,{\rm d}t\pm\int_{\Omega_T}\alpha\, p(\partial_t\sigma)(\triangle\sigma)\, z^2\,{\rm d}x\,{\rm d}t\\ \nonumber
 &-\int_{\Omega_T}\left[\frac{1}{2}|\nabla\sigma|^2\nabla(p^2)\cdot\nabla\sigma+p^2\nabla(|\nabla\sigma|^2)\cdot\nabla\sigma\right]\, z^2\,{\rm d}x\,{\rm d}t\\ 
 &+\int_{\Omega_T}\left[\alpha\,p^2|\nabla\sigma|^2\triangle\sigma-\alpha^2p^2|\triangle\sigma|^2)\right] z^2\,{\rm d}x\,{\rm d}t,
 \end{align}
 where $\partial_\nu=\vec{\nu}\cdot\nabla$, $\vec\nu=\vec\nu(x)=(\nu_1(x),\ldots,\nu_n(x))$ being the external unit normal vector on the boundary $\partial\Omega$ at $x$, and $dS$ denotes the trace of the Lebesgue measure on $\partial\Omega$.
  \end{lemma}
   For the proof, see Subsection 3.2.

 {\it 3.1.b. First Step: Lower  Bound of the Quantity I.}

 We will decompose the term $I$ in (\ref{I}) into several terms that we handle separately.
  $$I=I_l+I_{nl}+I_t.$$
  Here, we define $I_l$ as the sum of the integrals linear in $\sigma$ which do not have any time derivative, $I_{nl}$ as the sum of the integrals involving non-linear terms in $\sigma$  and without any
time derivative, and $I_t$ as the terms involving the time derivatives in $\sigma$. $I_l$, $I_{nl}$, and $I_t$ are as in 
(\ref{I,l}), ( \ref{I,nl}), and (\ref{I,t}).

About
 \begin{align}\label{I,l}
 \nonumber I_l&=\int_{\Sigma_T}|\partial_\nu z|^2p^2(\nabla\sigma\cdot\vec\nu){\rm d}S\,{\rm d}t
 +\int_{\Omega_T}\left[p\,(\nabla p\cdot\nabla\sigma)-\alpha \,p^2\triangle\sigma\right]|\nabla z|^2\,{\rm d}x\,{\rm d}t\\ \nonumber
 &-2\int_{\Omega_T}p(\nabla z\cdot\nabla\sigma)(\nabla z\cdot\nabla p)\,{\rm d}x\,{\rm d}t-2\int_{\Omega_T}p^2\sum_{i,j=1}^n\partial_i\partial_j\sigma\partial_iz\partial_jz\,{\rm d}x\,{\rm d}t\\ \nonumber
 &+\int_{\Omega_T}\left[\frac{1}{4}\bigg(3\nabla (p^2)\cdot\nabla\alpha+(1+\alpha)\triangle(p^2)\bigg)+\frac{1}{2}p^2\triangle\alpha\right]\triangle\sigma\,z^2\,{\rm d}x\,{\rm d}t\\ \nonumber
 &+\int_{\Omega_T}\left(\frac{3}{4}(1+\alpha)\nabla(p^2)+p^2\nabla\alpha\right)\cdot\nabla(\triangle\sigma)\,z^2\,{\rm d}x\,{\rm d}t\\ \nonumber
 &+\frac{1}{2}\int_{\Omega_T}\left[(1+\alpha)p^2\triangle(\triangle\sigma)+ \nabla p\cdot\nabla(\nabla p\cdot\nabla\sigma)+p\,\triangle(\nabla p\cdot\nabla\sigma)\right]\, z^2\,{\rm d}x\,{\rm d}t\\ 
 &+2\int_{\Omega_T}p\,(\nabla\sigma\cdot x)\frac{\mu}{|x|^4}\,\, z^2\,{\rm d}x\,{\rm d}t
 +\int_{\Omega_T}\alpha\, p\,\triangle\sigma\frac{\mu}{|x|^2}\,\, z^2\,{\rm d}x\,{\rm d}t,
 \end{align}
  we have the following estimate:
\begin{lemma}
There exist some positive constants such that for $\lambda>0$ large enough, we have:
\begin{align}\label{lem I,l}
\nonumber I_l\geq&\,\, \frac{C_0}{2}s\lambda^2\int_{\widetilde{\Theta}_T}\theta\phi|\nabla z|^2\,{\rm d}x\,{\rm d}t+C_1s\int_{\Omega_T}\theta|\nabla z|^2\,{\rm d}x\,{\rm d}t
 +C_2s\int_{\Omega_T}\theta\frac{|z|^2}{|x|^\gamma}\,{\rm d}x\,{\rm d}t\\ 
&-C_4s\lambda^2\int_{\omega_{0,T}}\theta\phi|\nabla z|^2\,{\rm d}x\,{\rm d}t
-C_5s\lambda^4\int_{\Omega_T}\theta \,z^2\,{\rm d}x\,{\rm d}t-C_6s\lambda^4\int_{\widetilde{\Theta}_T}\theta\phi\, z^2\,{\rm d}x\,{\rm d}t,
\end{align}
where $C_0=\frac{1}{2n}p_1^2\delta^2>0$, $C_1=2p_1^2-3p_1p_3r-2p_2\frac{\mu}{\mu^*}>0$, and $C_2=2p_2\frac{\mu}{\mu^*}>0$. 
\end{lemma}

The proof is given in Subsection 3.2. Note that the proof of Lemma 3.2 uses an improved form of the Hardy inequality (\ref{1.3}), which can be found for instance in \cite{M}, namely:
\begin{lemma}
For all $l\geq0$ and for all $0<\gamma<2$, there exists a positive constant $K_0=K_0(\gamma,l)>0$ such that
\begin{gather}\label{IH}
\mu^*\int_{\Omega}\frac{|z|^2}{|x|^2}\,{\rm d}x+l\int_{\Omega}\frac{|z|^2}{|x|^\gamma}\,{\rm d}x\leq\int_{\Omega}|\nabla z|^2\,{\rm d}x+K_0\int_{\Omega}|z|^2\,{\rm d}x,\qquad \forall z\in H_0^1(\Omega).
\end{gather}
\end{lemma}
We then consider $I_{nl}$.
\begin{align}\label{I,nl}
 \nonumber I_{nl}&=-\int_{\Omega_T}\left[\frac{1}{2}|\nabla\sigma|^2\nabla(p^2)\cdot\nabla\sigma+p^2\nabla(|\nabla\sigma|^2)\cdot\nabla\sigma\right]\, z^2\,{\rm d}x\,{\rm d}t\\ 
 &+\int_{\Omega_T}\left[\alpha\,p^2|\nabla\sigma|^2\triangle\sigma-\alpha^2p^2|\triangle\sigma|^2)\right] z^2\,{\rm d}x\,{\rm d}t.
 \end{align}
 Using (\ref{sigma}), we shall prove in Subsection 3.2 that
\begin{lemma}
For any $\lambda>0$ large enough, there exists  $s_1(\lambda)>0$ such that for $s \geq s_1(\lambda)$,
\begin{align}\label{lem I,nl}
\nonumber I_{nl}\geq&\, \,\frac{C_7}{2}s^3\lambda^4\int_{\widetilde{\Theta}_T}\theta^3\phi^3|z|^2\,{\rm d}x\,{\rm d}t+C_8s^3\int_{\Omega_T}\theta^3|x|^2|z|^2\,{\rm d}x\,{\rm d}t\\ 
&-C_9s^3\lambda^4\int_{\omega_{0,T}}\theta^3\phi^3|z|^2\,{\rm d}x\,{\rm d}t.
\end{align}
\end{lemma}
We finally estimate $I_t$.
 \begin{align}\label{I,t}
 I_{t}=&-\frac{1}{2}\int_{\Omega_T}\partial_t^2\sigma z^2\,{\rm d}x\,{\rm d}t
 \mp2\int_{\Omega_T} p\,(\nabla\sigma\cdot\nabla\partial_t\sigma)\, z^2\,{\rm d}x\,{\rm d}t\pm\int_{\Omega_T}\alpha\, p\,\partial_t\sigma\triangle\sigma\, z^2\,{\rm d}x\,{\rm d}t.
 \end{align}
 We further add to $I_t$ the last two integrals appearing in (\ref{lem I,l}) that we want to get rid of and define
\begin{gather}\label{I,r}
 I_r=I_t-C_5s\lambda^4\int_{\Omega_T}\theta \,z^2\,{\rm d}x\,{\rm d}t-C_6s\lambda^4\int_{\widetilde{\Theta}_T}\theta\phi \,z^2\,{\rm d}x\,{\rm d}t.
 \end{gather}
 In Subsection  3.2, we  shall prove that 
 \begin{lemma}
 For any $\lambda>0$ large enough, there exists $s_0(\lambda)\geq s_1(\lambda)$ such that for $s\geq s_0(\lambda)$,
 \begin{align}\label{lem I,r}
\nonumber |I_r|\leq&\,\,\frac{C_2}{2}s\int_{\Omega_T}\theta\frac{z^2}{|x|^\gamma}\,{\rm d}x\,{\rm d}t
+\frac{C_7}{4}s^3\lambda^4\int_{\widetilde{\Theta}_T}\theta^3\phi^3\,z^2\,{\rm d}x\,{\rm d}t\\ 
&+\frac{C_8}{2}s^3\int_{\Omega_T}\theta^3|x|^2\,z^2\,{\rm d}x\,{\rm d}t,
\end{align}
where $C_2$ is as in (\ref{lem I,l}), and $C_7$, $C_8$ are as in (\ref{lem I,nl}).
\end{lemma}
Noting (\ref{I,r}) and applying   (\ref{lem I,l}), (\ref{lem I,nl}) and (\ref{lem I,r}) in (\ref{I}), we obtain a lower bound of $I$:
\begin{align}\label{II}
\nonumber I&\geq\,\, \frac{C_0}{2}s\lambda^2\int_{\widetilde{\Theta}_T}\theta\phi|\nabla z|^2\,{\rm d}x\,{\rm d}t+C_1s\int_{\Omega_T}\theta|\nabla z|^2\,{\rm d}x\,{\rm d}t+\frac{C_2}{2}s\int_{\Omega_T}\theta\frac{|z|^2}{|x|^\gamma}\,{\rm d}x\,{\rm d}t\\ \nonumber
&+\frac{C_7}{4}s^3\lambda^4\int_{\widetilde{\Theta}_T}\theta^3\phi^3|z|^2\,{\rm d}x\,{\rm d}t+\frac{C_8}{2}s^3\int_{\Omega_T}\theta^3|x|^2|z|^2\,{\rm d}x\,{\rm d}t\\ 
&-C_9s^3\lambda^4\int_{\omega_{0,T}}\theta^3\phi^3|z|^2\,{\rm d}x\,{\rm d}t-C_4s\lambda^2\int_{\omega_{0,T}}\theta\phi|\nabla z|^2\,{\rm d}x\,{\rm d}t.
\end{align}

{\it 3.1.c. Second Step: Estimate of the Time Derivative.}

Next, we proceed to the estimate of the time derivative of the solution. We recall that, by (\ref{2.13}), we have
\begin{gather}\label{V.6}
\pm\partial_tz=Az-2p\,(\nabla\sigma\cdot\nabla z)-\left(1+\alpha\right)p (\triangle\sigma)z
-(\nabla p\cdot\nabla\sigma)z
\end{gather}
In order to estimate $z_t$, we proceed separately in $B_r\times(0,T)$ and in $\widetilde{\Theta}_T$.

Firstly,  we note that, $\forall (x,t)\in\widetilde{\Theta}_T$, $j, k=1,...,n$,
\begin{align}\label{2.27}
 |\partial_j\sigma|\leq Cs\lambda\theta\phi,\quad |\partial_k\partial_j\sigma|\leq Cs\lambda^2\theta\phi,\quad
|\nabla\triangle\sigma|\leq Cs\lambda^3\theta\phi,\quad |\triangle^2\sigma|\leq Cs\lambda^4\theta\phi.
\end{align}
Here and henceforth, all the constants $C$ are positive and in dependent of $s$ and $\lambda$. We will verify (\ref{2.27}) in the appendix. For $s\geq1$ and $\lambda>0$ large enough but independent of $s$,  (\ref{V.6}) and (\ref{2.27}) yield
\begin{align}\label{V.10}
\nonumber \frac{1}{s\lambda}\int_{\widetilde{\Theta}_T}&\frac{1}{\theta\phi}z_t^2\,{\rm d}x\,{\rm d}t\leq \frac{C}{s\lambda}\int_{\widetilde{\Theta}_T}\frac{1}{\theta\phi}\vert{Az}\vert^2\,{\rm d}x\,{\rm d}t
+Cs\lambda\int_{\widetilde{\Theta}_T}\theta\phi|\nabla z|^2\,{\rm d}x\,{\rm d}t\\ \nonumber
&+Cs\lambda^3\int_{\widetilde{\Theta}_T}\theta\phi z^2\,{\rm d}x\,{\rm d}t+Cs\lambda\int_{\widetilde{\Theta}_T}\theta\phi z^2\,{\rm d}x\,{\rm d}t\\ 
&\leq \frac{1}{4}\Vert{Az}\Vert^2
+\frac{C_0}{4}s\lambda^2\int_{\widetilde{\Theta}_T}\theta\phi|\nabla z|^2\,{\rm d}x\,{\rm d}t
+\frac{C_7}{8}s^3\lambda^4\int_{\widetilde{\Theta}_T}\theta^3\phi^3 z^2\,{\rm d}x\,{\rm d}t,
\end{align}
where $C_0$ and $C_8$ are  constants appearing in (\ref{II}), and $\Vert{\cdot}\Vert$ denotes as in Subsection 3.1.a.

In $B_{r}\times(0,T)$, all the  computations are  explicit. For $\forall \lambda\geq2 $, by Appendix, we have
\begin{align}\label{2.28}
|\partial_j\sigma|\leq Cs\lambda\theta|x|,\quad |\partial_k\partial_j\sigma|\leq Cs\lambda^2\theta,
\quad |\nabla\triangle\sigma|\leq Cs\lambda^3\theta,\quad |\triangle^2\sigma|\leq Cs\lambda^4\theta.
\end{align}
Then, by (\ref{V.6}) and (\ref{2.28}), we deduce that
\begin{align*}
\frac{1}{s\lambda^5}\int_{B_r\times(0,T)}&\frac{1}{\theta}z_t^2\,{\rm d}x\,{\rm d}t
\leq \frac{C}{s\lambda^5}\int_{B_r\times(0,T)}\frac{1}{\theta}|Az|^2\,{\rm d}x\,{\rm d}t
+\frac{C}{\lambda^3}s\int_{B_r\times(0,T)}\theta |x|^2|\nabla z|^2\,{\rm d}x\,{\rm d}t\\ \nonumber
&+\frac{C}{\lambda}s\int_{B_r\times(0,T)}\theta z^2\,{\rm d}x\,{\rm d}t
+\frac{C}{\lambda^3}s\int_{B_r\times(0,T)}\theta |x|^2 z^2\,{\rm d}x\,{\rm d}t\\ \nonumber
&\leq \frac{C}{s\lambda^5}\Vert{Az}\Vert^2+\frac{Cr^2}{\lambda^3}s\int_{\Omega_T}\theta|\nabla z|^2\,{\rm d}x\,{\rm d}t+\frac{C}{\lambda}s\int_{\Omega_T}\theta\frac{z^2}{|x|^\gamma}\,{\rm d}x\,{\rm d}t.
\end{align*}
Finally, we get for $s\geq1$ and $\lambda>0$ large enough but independent of $s$:
\begin{align}\label{V.8}
 \frac{1}{s\lambda^5}\int_{B_r\times(0,T)}\frac{1}{\theta}z_t^2\,{\rm d}x\,{\rm d}t\leq &
\frac{1}{4}\Vert{Az}\Vert^2+\frac{C_1}{2}s\int_{\Omega_T}\theta|\nabla z|^2\,{\rm d}x\,{\rm d}t+\frac{C_2}{4}s\int_{\Omega_T}\theta\frac{z^2}{|x|^\gamma}\,{\rm d}x\,{\rm d}t,
\end{align}
where $C_2$ is the constant appearing in (\ref{II}).

{\it 3.1.d. Partial Conclusion: Estimates in z.}

By (\ref{2.14}), (\ref{I.0}), (\ref{II}), (\ref{V.10}) and (\ref{V.8}),  we conclude that
\begin{align}\label{V.11}
\nonumber &\frac{C_0}{4}s\lambda^2\int_{\widetilde{\Theta}_T}\theta\phi|\nabla z|^2\,{\rm d}x\,{\rm d}t
+\frac{C_2}{4}s\int_{\Omega_T}\theta\frac{|z|^2}{|x|^\gamma}\,{\rm d}x\,{\rm d}t+\frac{C_1}{2}s\int_{\Omega_T}\theta|\nabla z|^2\,{\rm d}x\,{\rm d}t\\ \nonumber
&+\frac{C_7}{8}s^3\lambda^4\int_{\widetilde{\Theta}_T}\theta^3\phi^3|z|^2\,{\rm d}x\,{\rm d}t+\frac{C_8}{2}s^3\int_{\Omega_T}\theta^3|x|^2|z|^2\,{\rm d}x\,{\rm d}t\\ \nonumber
&+\frac{1}{s\lambda^5}\int_{B_r\times(0,T)}\frac{1}{\theta}z_t^2\,{\rm d}x\,{\rm d}t
+\frac{1}{s\lambda}\int_{\widetilde{\Theta}_T}\frac{1}{\theta\phi}z_t^2\,{\rm d}x\,{\rm d}t\\ \nonumber
&\leq I+\frac{1}{2}\Vert{Az}\Vert^2+C_4s\lambda^2\int_{\omega_{0,T}}\theta\phi|\nabla z|^2\,{\rm d}x\,{\rm d}t+C_9s^3\lambda^4\int_{\omega_{0,T}}\theta^3\phi^3|z|^2\,{\rm d}x\,{\rm d}t\\ 
&\leq \Vert{{\rm e}^{-\sigma}g}\Vert^2+C_4s\lambda^2\int_{\omega_{0,T}}\theta\phi|\nabla z|^2\,{\rm d}x\,{\rm d}t+C_9s^3\lambda^4\int_{\omega_{0,T}}\theta^3\phi^3|z|^2\,{\rm d}x\,{\rm d}t,
\end{align}
where we used (\ref{2.14}) in the last inequality.

{\it 3.1.e. Final Conclusion: Estimates in w.}

To conclude the proof of Theorem 2.1, it now remains to undo the change of variable $z=e^{-\sigma}w$ and deduce from (\ref{V.11}) the required estimate in terms of $w$. In fact, the proof is similar to Subsection V.1.e in \cite{J} except the following Lemma 3.7.

First of all, we compute
\begin{gather}\label{V.12}
{\rm e}^{-2\sigma}|\nabla w|^2=|\nabla z|^2+|\nabla\sigma|^2z^2+2(\nabla \sigma\cdot\nabla z)z\leq 2|\nabla z|^2+2|\nabla\sigma|^2z^2.
\end{gather}
Thus, by (\ref{2.27}) and (\ref{2.28}), we have
\begin{gather}\label{3.27}
s\lambda^2\int_{\widetilde{\Theta}_T}\theta\phi {\rm e}^{-2\sigma}|\nabla w|^2\,{\rm d}x\,{\rm d}t\leq Cs^3\lambda^4\int_{\widetilde{\Theta}_T}\theta^3\phi^3|z|^2\,{\rm d}x\,{\rm d}t
+Cs\lambda^2\int_{\widetilde{\Theta}_T}\theta\phi|\nabla z|^2\,{\rm d}x\,{\rm d}t,
\end{gather}
\begin{gather}\label{3.28}
\frac{s}{\lambda^2}\int_{B_r\times(0,T)}\theta {\rm e}^{-2\sigma}|\nabla w|^2\,{\rm d}x\,{\rm d}t\leq Cs^3\int_{B_r\times(0,T)}\theta^3|x|^2|z|^2\,{\rm d}x\,{\rm d}t
+C\frac{s}{\lambda^2}\int_{B_r\times(0,T)}\theta|\nabla z|^2\,{\rm d}x\,{\rm d}t,
\end{gather}
which can be estimates by the right hand side of (\ref{V.11}).
Next we compute
\begin{gather*}
{\rm e}^{-\sigma}w_t=z_t+s\theta'({\rm e}^{2\lambda\sup\psi}-\frac{1}{2}|x|^2-{\rm e}^{\lambda\psi})z.
\end{gather*}
Similarly to (V.19)-(V.21) in \cite{J}, using the fact that $|\theta'|\leq C\theta^{1+1/{k_0}}$ and Young inequality, taking $k_0=1+2/\gamma$ and $\lambda>0$ large enough but independent of s , we can obtain 
\begin{align}\label{107}
\nonumber \frac{1}{s}&{\rm e}^{-4\lambda\sup\psi}\int_{B_r\times(0,T)}\frac{1}{\theta}{\rm e}^{-2\sigma}w_t^2\,{\rm d}x\,{\rm d}t\\ 
&\leq \frac{C}{s}{\rm e}^{-4\lambda\sup\psi}\int_{B_r\times(0,T)}\frac{1}{\theta}z_t^2\,{\rm d}x\,{\rm d}t
+Cs\int_{\Omega_T}\theta^3|x|^{2}z^{2}\,{\rm d}x\,{\rm d}t+Cs\int_{\Omega_T}\theta\frac{|z|^2}{|x|^\gamma}\,{\rm d}x\,{\rm d}t,
\end{align}
which can be estimated by the right hand side of (\ref{V.11}).

For $(x,t)\in\widetilde{\Theta}$, similarly to (V.18) in \cite{J}, by  $\phi\geq1$, and $k_0>1,$ we have
\begin{align}\label{1077}
\nonumber \frac{1}{s}{\rm e}^{-4\lambda\sup\psi}\int_{\widetilde{\Theta}_T}\frac{1}{\theta\phi}{\rm e}^{-2\sigma}w_t^2\,{\rm d}x\,{\rm d}t&\leq \frac{C}{s}{\rm e}^{-4\lambda\sup\psi}\int_{\widetilde{\Theta}_T}\frac{1}{\theta\phi}z_t^2\,{\rm d}x\,{\rm d}t
+Cs\int_{\widetilde{\Theta}_T}\frac{1}{\phi}\theta^{1+2/{k_0}}z^2\,{\rm d}x\,{\rm d}t,\\
&\leq\frac{C}{s}{\rm e}^{-4\lambda\sup\psi}\int_{\widetilde{\Theta}_T}\frac{1}{\theta\phi}z_t^2\,{\rm d}x\,{\rm d}t
+Cs\int_{\widetilde{\Theta}_T}\theta^3\phi^3z^2\,{\rm d}x\,{\rm d}t,
\end{align}
which can be estimates by the right hand side of (\ref{V.11}).

On the other hand, by (\ref{V.12}), we have
\begin{gather*}
|\nabla z|^2= {\rm e}^{-2\sigma}|\nabla w|^2-|\nabla\sigma|^2z^2-2(\nabla\sigma\cdot\nabla z)\,z .
\end{gather*}
Hence
\begin{gather*}
|\nabla z|^2\leq {\rm e}^{-2\sigma}|\nabla w|^2+|\nabla\sigma|^2z^2+\frac{1}{2}|\nabla z|^2 +C|\nabla\sigma|^2z^2.
\end{gather*}
Hence, using (\ref{2.27}), we get
\begin{gather*}
\forall (x,t)\in\omega_{0,T},\quad \frac{1}{2} |\nabla z|^2\leq {\rm e}^{-2\sigma}|\nabla w|^2+Cs^2\lambda^2\theta^2\phi^2z^2.
\end{gather*}
Thus,
\begin{gather}\label{10777}
s\lambda^2\int_{\omega_{0,T}}\theta\phi |\nabla z|^2\,{\rm d}x\,{\rm d}t\leq s\lambda^2\int_{\omega_{0,T}}\theta\phi {\rm e}^{-2\sigma}|\nabla w|^2\,{\rm d}x\,{\rm d}t+Cs^3\lambda^4\int_{\omega_{0,T}}\theta^3\phi^3z^2\,{\rm d}x\,{\rm d}t.
\end{gather}
Finally, by using the fact that $\omega_{0}\subset\omega$ and the fact that $z={\rm e}^{-\sigma}w$, the combination of (\ref{3.27})- -(\ref{10777}) yields
\begin{align}\label{3.53}
\nonumber&s\lambda^2\int_{\widetilde{\Theta}_T}\theta\phi {\rm e}^{-2\sigma}|\nabla w|^2\,{\rm d}x\,{\rm d}t
+\frac{s}{\lambda^2}\int_{\Omega_T}\theta {\rm e}^{-2\sigma}|\nabla w|^2\,{\rm d}x\,{\rm d}t\\ \nonumber
&+s\int_{\Omega_T}\theta {\rm e}^{-2\sigma}\frac{|w|^2}{|x|^\gamma}\,{\rm d}x\,{\rm d}t+s^3\lambda^4\int_{\widetilde{\Theta}_T}\theta^3\phi^3{\rm e}^{-2\sigma}|w|^2\,{\rm d}x\,{\rm d}t
+s^3\int_{\Omega_T}\theta^3{\rm e}^{-2\sigma}|x|^2|w|^2\,{\rm d}x\,{\rm d}t\\ \nonumber
&+\frac{1}{s}{\rm e}^{-4\lambda\sup\psi}\int_{B_r\times(0,T)}\frac{1}{\theta}{\rm e}^{-2\sigma}w_t^2\,{\rm d}x\,{\rm d}t
+\frac{1}{s}{\rm e}^{-4\lambda\sup\psi}\int_{\widetilde{\Theta}_T}\frac{1}{\theta\phi}{\rm e}^{-2\sigma}w_t^2\,{\rm d}x\,{\rm d}t\\ \nonumber
&\leq C\int_{\Omega_T}g^2{\rm e}^{-2\sigma}\,{\rm d}x\,{\rm d}t+Cs\lambda^2\int_{\omega_{0,T}}\theta\phi {\rm e}^{-2\sigma}|\nabla w|^2\,{\rm d}x\,{\rm d}t\\ 
&+Cs^3\lambda^4\int_{\omega_{T}}\theta^3\phi^3{\rm e}^{-2\sigma}|w|^2\,{\rm d}x\,{\rm d}t.
\end{align}
Next, we use the following Lemma Inequality:
\begin{lemma}\rm{(Cacciopoli's Inequality).}
For $\sigma, \psi,\theta,k_0$ and $\phi$ be defined as above in this paper,
there exists a constant $C>0$ $($independent of $0\leq\mu\leq p_1\mu^*(n)$ and of $s$ and $\lambda$$)$ such that  for any  $\lambda>0$ large enough and for any $s>1$, any solution w of (\ref{pm}) satisfies
\begin{align}\label{2.26}
\int_{\omega_{0,T}}\theta\phi {\rm e}^{-2{\sigma}}|\nabla w|^2\,{\rm d}x{\rm d}t
\leq\, Cs^2{\rm e}^{2\lambda\sup{\psi}}\int_{\omega_{T}}\theta^3\phi^3{\rm e}^{-2{\sigma}}w^2\,{\rm d}x{\rm d}t+\frac{1}{s\lambda^2}\int_{\omega_{T}}{\rm e}^{-2{\sigma}}g^2\,{\rm d}x{\rm d}t.
\end{align}
\end{lemma}
The proof of lemma 3.7 is similar to the proof of  Lemma III.3  in \cite{VZ1}. For completeness, we will prove Lemma 3.7  in Subsection 3.2. Thus, using (\ref{3.53}) and (\ref{2.26}), we deduce (\ref{CE}).

This ends the proof of Theorem 2.1.

\subsection{Proofs of Technical Lemmas}

Here we present the proofs of the technical Lemmas stated in Subsection 3.1. 

{\it Proof of Lemma 3.1.}

To make the computations easier, we define
\begin{gather*}
S_1z={\rm div}(p\nabla z), \quad S_2z=p\,z|\nabla\sigma|^2,
\quad S_3z=\pm\,\sigma_t\,z,\quad S_4z=\frac{\mu}{|x|^2}\,z,\\ \nonumber
 A_1z=\pm\,\partial_tz,\quad A_2z=2p\,(\nabla\sigma\cdot\nabla z),\quad
A_3z=\left(1+\alpha\right)p(\triangle\sigma)z,\quad A_4z=(\nabla p\cdot\nabla\sigma)z,
\end{gather*}
and denote by $I_{ij}$ the scalar product $\left<S_i,A_j\right>$. We will compute each term using integration by parts and the boundary conditions (\ref{boundary1}) and (\ref{boundary2}). Thus,
 $$I=\sum_{k,j=1}^4I_{kj}-\int_{\Omega_T}\alpha^2p^2(\triangle\sigma)^2z^2\,{\rm d}x\,{\rm d}t. $$
Computation of $I_{11}$: By (\ref{boundary1}) and (\ref{boundary2}), we have
\begin{align*}
I_{11}=\pm\left(-\frac{1}{2}\int_{\Omega_T}\partial_t(p\,|\nabla z|^2)\,{\rm d}x\,{\rm d}t\right)=0.
\end{align*}\\
Computation of $I_{12}$:
Note that, since $z$ vanishes on the boundary, its gradient $\nabla z$ on the boundary is normal to the boundary, and therefore $\nabla z=\partial_\nu z \,\, \vec{\nu}$ on $\Sigma_T$, where $\vec{\nu}$ denotes the unit outward normal vector on the boundary.
\begin{align*}
&I_{12}=\int_{\Sigma_T}p^2(\nabla\sigma\cdot\vec\nu)|\partial_\nu z|^2\,{\rm d}S\,{\rm d}t
-2\int_{\Omega_T}p\,(\nabla\sigma\cdot\nabla z)(\nabla p\cdot\nabla z)\,{\rm d}x\,{\rm d}t\\
&+2\int_{\Omega_T}p\,(\nabla\sigma\cdot\nabla p)|\nabla z|^2\,{\rm d}x\,{\rm d}t
+\int_{\Omega_T}p^2(\triangle\sigma)|\nabla z|^2\,{\rm d}x\,{\rm d}t-2\int_{\Omega_T}p^2\sum_{k,j=1}^{n}(\partial_k\partial_j\sigma)
(\partial_kz)\partial_jz\,{\rm d}x\,{\rm d}t.
\end{align*}
Computation of $I_{13}$: Integrating by parts and using (\ref{boundary2}), we have
\begin{align*}
&I_{13}=\int_{\Omega_T}\left[\frac{1}{4}\bigg(3\nabla (p^2)\cdot\nabla\alpha+(1+\alpha)\triangle(p^2)\bigg)+\frac{1}{2}p^2\triangle\alpha\right](\triangle\sigma)\,z^2\,{\rm d}x\,{\rm d}t\\ \nonumber
 &+\int_{\Omega_T}\left(\frac{3}{4}(1+\alpha)\nabla(p^2)\cdot\nabla(\triangle\sigma)+p^2\nabla\alpha\cdot\nabla(\triangle\sigma)\right)\,z^2\,{\rm d}x\,{\rm d}t\\
 &+\frac{1}{2}\int_{\Omega_T}(1+\alpha)p^2(\triangle^2\sigma)\, z^2\,{\rm d}x\,{\rm d}t-\int_{\Omega_T}(1+\alpha)p^2(\triangle\sigma)\,|\nabla z|^2\,{\rm d}x\,{\rm d}t.
\end{align*}
Computation of $I_{14}$:
\begin{align*}
&I_{14}=\frac{1}{2}\int_{\Omega_T}\left[\nabla p\cdot\nabla(\nabla p\cdot\nabla\sigma)+p\,\triangle(\nabla p\cdot\nabla\sigma)\right]\, z^2\,{\rm d}x\,{\rm d}t
-\int_{\Omega_T}p\,(\nabla p\cdot\nabla\sigma)|\nabla z|^2\,{\rm d}x\,{\rm d}t.
 \end{align*}
Computation of $I_{21}$:
\begin{align*}
I_{21}&=\mp\int_{\Omega_T}p\,(\nabla\sigma\cdot\nabla\partial_t\sigma) z^2\,{\rm d}x\,{\rm d}t.
\end{align*}
Computation of $I_{22}$:
\begin{align*}
&I_{22}=-\int_{\Omega_T} \left[|\nabla\sigma|^2\nabla(p^2)\cdot\nabla\sigma+p^2\nabla(|\nabla\sigma|^2)\cdot\nabla\sigma\right]z^2\,{\rm d}x\,{\rm d}t
 -\int_{\Omega_T}\,p^2|\nabla\sigma|^2(\triangle\sigma)z^2 \,{\rm d}x\,{\rm d}t.
\end{align*}
Computation of $I_{23}$:
\begin{align*}
I_{23}=\int_{\Omega_T}(1+\alpha)p^2|\nabla\sigma|^2(\triangle\sigma)\, z^2\,{\rm d}x\,{\rm d}t.
\end{align*}
Computation of $I_{24}$:
\begin{align*}
I_{24}=\int_{\Omega_T}p(\nabla p\cdot\nabla\sigma)|\nabla\sigma|^2 z^2\,{\rm d}x\,{\rm d}t.
\end{align*}
Computation of $I_{31}$:
\begin{align*}
I_{31}&=-\frac{1}{2}\int_{\Omega_T}\left(\partial_{t}^2\sigma\right)\, z^2\,{\rm d}x\,{\rm d}t.
\end{align*}
Computation of $I_{32}$: Integrating by parts and using (\ref{boundary2}), we have
\begin{align*}
I_{32}&=\mp\int_{\Omega_T}\left[\left(\nabla p\cdot\nabla\sigma\right)\partial_t\sigma\
+p\,\left(\triangle \sigma\right)\partial_t\sigma+p\,\left(\nabla\sigma\cdot\nabla\partial_t\sigma\right)\right]  z^2\,{\rm d}x\,{\rm d}t.
\end{align*}
Computation of $I_{33}$:
\begin{align*}
I_{33}=\pm\int_{\Omega_T}(1+\alpha)p(\partial_t\sigma)\left(\triangle\sigma\right) z^2\,{\rm d}x\,{\rm d}t.
\end{align*}
Computation of $I_{34}$:
\begin{gather*}
I_{34}=\pm\int_{\Omega_T}(\partial_t\sigma)\left(\nabla\sigma\cdot\nabla p\right) z^2\,{\rm d}x\,{\rm d}t.
\end{gather*}
Computation of $I_{41}$:
\begin{gather*}
I_{41}=\pm\int_{\Omega_T}\frac{1}{2}\partial_t\left(\frac{\mu}{|x|^2}\,z^2\right) \,{\rm d}x\,{\rm d}t=0.
\end{gather*}
Computation of $I_{42}$:
\begin{align*}
I_{42}&=\int_{\Omega_T}p\,\left(\nabla\sigma\cdot x \right)\frac{2\mu}{|x|^4}z^2\,{\rm d}x\,{\rm d}t-\int_{\Omega_T}\left(\nabla\sigma\cdot\nabla p\right)\frac{\mu}{|x|^2}z^2\,{\rm d}x\,{\rm d}t
-\int_{\Omega_T}p\,\left(\triangle\sigma\right)\frac{\mu}{|x|^2}z^2\,{\rm d}x\,{\rm d}t.\\ \nonumber
\end{align*}
Computation of $I_{43}$:
\begin{align*}
I_{43}=\int_{\Omega_T}(1+\alpha)p\,\left(\triangle\sigma\right) \frac{\mu}{|x|^2}\,z^2\,{\rm d}x\,{\rm d}t.
\end{align*}
Computation of $I_{44}$:
\begin{align*}
I_{44}=\int_{\Omega_T}\left(\nabla\sigma\cdot\nabla p\right)\frac{\mu}{|x|^2}\,z^2\,{\rm d}x\,{\rm d}t.
\end{align*}

By these computations, we obtain (\ref{I}) in Lemma 3.1.\\

{\it Proof of Lemma 3.2.}

We decompose $I_l$ as
\begin{equation*}
I_l=I_{l,1}+I_{l,2},
\end{equation*}
where $I_{l,1}$ and $I_{l,2}$ are the terms in $I_l$ corresponding respectively in $\widetilde{\Theta}_T$ and in $B_r\times(0,T)$.

First, let us compute  $I_{l,1}$. By the appendix, we have 
\begin{align}\label{I,l1,1}
 \nonumber I_{l,1}&\geq\,-s\int_{\Sigma_T}\theta\,p^2\left[x\cdot\nu+\lambda\phi(\nabla\psi\cdot\nu)\right]|\partial_\nu z|^2\,{\rm d}S\,{\rm d}t\\ \nonumber
&+s\int_{\widetilde{\Theta}_T}\alpha\,p^2\theta(n+\lambda\phi\triangle\psi+\lambda^2\phi |\nabla\psi|^2)|\nabla z|^2\,{\rm d}x\,{\rm d}t+2s\lambda^2\int_{\widetilde{\Theta}_T} \phi\,p^2\theta(\nabla\psi\cdot\nabla z)^2\,{\rm d}x\,{\rm d}t\\ \nonumber
&+2s\int_{\widetilde{\Theta}_T}p\,\theta[x\cdot\nabla z+\lambda\phi(\nabla\psi\cdot\nabla z)](\nabla p\cdot\nabla z)\,{\rm d}x\,{\rm d}t\\ \nonumber
&-s\int_{\widetilde{\Theta}_T}p\,\theta[x\cdot\nabla p+\lambda\phi(\nabla\psi\cdot\nabla p)]|\nabla z|^2\,{\rm d}x\,{\rm d}t\\
&+2s\int_{\widetilde{\Theta}_T}p^2\theta\sum_{i,j=1}^n(\delta_{ij}+\lambda\phi\partial_i\partial_j\psi)\partial_iz\partial_jz\,{\rm d}x\,{\rm d}t
-Cs\lambda^4\int_{\widetilde{\Theta}_T}\theta\,\phi\,| z|^2\,{\rm d}x\,{\rm d}t.
\end{align}
Besides, due to the particular choice of $\psi$, and especially (\ref{bounded}), the boundedness of $\widetilde{\Theta}$, and (\ref{a}),  (\ref{alpha}),  one can obtain following estimates:
\begin{gather}\label{es}
 \nonumber\alpha\,p^2|\nabla\psi|^2|\nabla z|^2\geq \frac{1}{n}\,p_1^2\delta^2|\nabla z|^2,\quad (x,t)\in\Theta_T,\\
p^2|\nabla\psi|^2|\nabla z|^2\leq C|\nabla z|^2,\quad (x,t)\in\omega_{0,T},\\ \nonumber
|\Delta^2\psi|+|\Delta\psi|+|\nabla\psi|+|\nabla\Delta\psi|\leq C,\quad (x,t)\in\widetilde{\Theta}_T.
\end{gather}
Then for $\lambda>0$ large enough, we have
\begin{align}\label{2.61}
I_{l,1}\geq& -s\int_{\Sigma_T}\theta\,p^2\left[x\cdot\nu+\lambda\phi(\nabla\psi\cdot\nu)\right]|\partial_\nu z|^2\,{\rm d}S\,{\rm d}t\\ \nonumber
+&C_0s\lambda^2\int_{\Theta_T}\theta\phi|\nabla z|^2\,{\rm d}x\,{\rm d}t
-Cs\lambda^4\int_{\widetilde{\Theta}_T}\theta\phi| z|^2\,{\rm d}x\,{\rm d}t-Cs\lambda^2\int_{\omega_{0,T}}\theta\phi|\nabla z|^2\,{\rm d}x\,{\rm d}t,
\end{align}
where $C_0=\frac{1}{2n}p_1^2\delta^2.$
Moreover, due to the properties (\ref{psi}) and (\ref{bounded})for $\lambda>0$ large enough,, the sum of boundary terms in the RHS of (\ref{2.61}) is positive. Indeed, from (\ref{psi}) and (\ref{bounded}), for $\forall (x,t)\in \Sigma_T$, $\nabla\psi\cdot\vec{\nu}=-|\nabla\psi|\leq-\delta$, $\phi=1$, and thus for $\lambda$ large enough,
\begin{align*}
x\cdot\nu+\lambda\phi(\nabla\psi\cdot\nu)&=x\cdot\nu-\lambda|\nabla\psi|
\leq x\cdot\nu-\lambda \delta\leq0.
\end{align*}
Thus,
\begin{align}\label{I,l1,2}
I_{l,1}\geq&\,C_0s\lambda^2\int_{\widetilde{\Theta}_T}\theta\phi|\nabla z|^2\,{\rm d}x\,{\rm d}t
-Cs\lambda^4\int_{\widetilde{\Theta}_T}\theta\phi| z|^2\,{\rm d}x\,{\rm d}t-Cs\lambda^2\int_{\omega_{0,T}}\theta\phi|\nabla z|^2\,{\rm d}x\,{\rm d}t.
\end{align}

Next, we compute $I_{l,2}$. In this case, using Appendix and the fact that $\alpha$ vanishes in $B(0,r/2)$ by (\ref{alpha}), for $\lambda\geq4$ large enough, all the computations are explicit:
\begin{align}\label{I,l2,1}
 \nonumber I_{l,2}&\geq\,2s\int_{B_r\times(0,T)}p\,\theta\left(1+\lambda\frac{|x|^{\lambda-2}}{r^\lambda}\right)(\nabla z\cdot\nabla p)(\nabla z\cdot x)\,{\rm d}x\,{\rm d}t\\ \nonumber
&+s\int_{B_r\times(0,T)}\alpha \,p^2\theta\left[n+(n-2)\lambda\frac{|x|^{\lambda-2}}{r^\lambda}+\lambda^2\frac{|x|^{\lambda-2}}{r^\lambda}\right]|\nabla z|^2\,{\rm d}x\,{\rm d}t\\ \nonumber
&-s\int_{B_r\times(0,T)}p\,\theta\left(1+\lambda\frac{|x|^{\lambda-2}}{r^\lambda}\right)(\nabla p\cdot x)|\nabla z|^2\,{\rm d}x\,{\rm d}t\\ \nonumber
&+2s\int_{B_r\times(0,T)}p^2\theta\left(1+\lambda\frac{|x|^{\lambda-2}}{r^\lambda}\right)|\nabla z|^2\,{\rm d}x\,{\rm d}t\\ \nonumber
&+2s\int_{B_r\times(0,T)}p^2\theta\,\lambda(\lambda-2)\frac{|x|^{\lambda-4}}{r^\lambda}|\nabla z\cdot x|^2\,{\rm d}x\,{\rm d}t\\
&-2s\int_{B_r\times(0,T)}p\,\theta\frac{\mu}{|x|^2}z^2\,{\rm d}x\,{\rm d}t-Cs\lambda^4\int_{B_r\times(0,T)}\theta\, z^2\,{\rm d}x\,{\rm d}t.
\end{align}

Noting that $r$ is a fixed number defined by  (\ref{r}), $\theta$ only depends on the time variable $t$ and  $\alpha$ is non-negative by $(\ref{alpha})$, and using (\ref{a}), we obtain that for $\lambda>2$ large enough, there exists a constant $C>0$ such that
\begin{align*}
\nonumber I_{l,2}&\geq s\int_{B_r\times(0,T)}\theta\left\{p\left(2p_1-3p_3 r\right)\left(1+\lambda\frac{|x|^{\lambda-2}}{r^\lambda}\right)|\nabla z|^2-2p_2\frac{\mu}{|x|^2}z^2\right\}\,{\rm d}x\,{\rm d}t\\
&-Cs\lambda^4\int_{B_r\times(0,T)}\theta z^2\,{\rm d}x\,{\rm d}t\\
&\geq s\int_{B_r\times(0,T)}\theta\left(p_1(2p_1-3p_3r)|\nabla z|^2
-2p_2\frac{\mu}{|x|^2}z^2\right)\,{\rm d}x\,{\rm d}t-Cs\lambda^4\int_{B_r\times(0,T)}\theta z^2\,{\rm d}x\,{\rm d}t\\
&\geq s\int_{\Omega_T}\theta\left(p_1(2p_1-3p_3r)|\nabla z|^2
-2p_2\frac{\mu}{|x|^2}z^2\right)\,{\rm d}x\,{\rm d}t\\ \nonumber
&-Cs\lambda^4\int_{\Omega_T}\theta z^2\,{\rm d}x\,{\rm d}t-Cs\int_{\widetilde{\Theta}_T}\theta|\nabla z|^2\,{\rm d}x\,{\rm d}t.
\end{align*}
Thus, from the Hardy improved inequality (\ref{IH}) with $l=1$, one obtain, for $\forall\,\, 0\leq\mu <\frac{p_1^2}{p_2}\mu^*<p_1\mu^*$,
\begin{gather*}
-2p_2\int_{\Omega}\mu\frac{z^2}{|x|^2}\,{\rm d}x\geq-2p_2\frac{\mu}{\mu^*}\int_{\Omega}|\nabla z|^2\,{\rm d}x+2p_2\frac{\mu}{\mu^*}\int_{\Omega}\frac{z^2}{|x|^\gamma}\,{\rm d}x-K_0p_2\frac{\mu}{\mu^*}\int_{\Omega}z^2\,{\rm d}x,
\end{gather*}
 Then for $\lambda>0$ large enough, one has
\begin{align}\label{I,l2,2}
\nonumber I_{l,2}&\geq s\int_{\Omega_T}\theta\left(2p_1^2-3p_1p_3r-2p_2\frac{\mu}{\mu^*}\right)|\nabla z|^2\,{\rm d}x\,{\rm d}t+2p_2\frac{\mu}{\mu^*}s\int_{\Omega_T}\theta\frac{z^2}{|x|^\gamma}\,{\rm d}x\,{\rm d}t\\ \nonumber
&-Cs\lambda^4\int_{\Omega_T}\theta z^2\,{\rm d}x\,{\rm d}t-Cs\int_{\widetilde{\Theta}_T}\theta|\nabla z|^2\,{\rm d}x\,{\rm d}t\\ \nonumber
&= C_1s\int_{\Omega_T}\theta|\nabla z|^2\,{\rm d}x\,{\rm d}t
+C_2s\int_{\Omega_T}\theta\frac{z^2}{|x|^\gamma}\,{\rm d}x\,{\rm d}t\\
&-Cs\lambda^4\int_{\Omega_T}\theta z^2\,{\rm d}x\,{\rm d}t-Cs\int_{\widetilde{\Theta}_T}\theta|\nabla z|^2\,{\rm d}x\,{\rm d}t.
\end{align}
Here, $C_1=2p_1^2-3p_1p_3r-2p_2\frac{\mu}{\mu^*}$, $C_2= 2p_2\frac{\mu}{\mu^*}>0$.

Hence, for $\lambda>0$ large enough, combining (\ref{I,l1,2}) and (\ref{I,l2,2}) gives Lemma 3.2.

{\it Proof of Lemma 3.4.}

Similarly, we denote $I_{nl}=I_{nl,1}+I_{nl,2}$, where $I_{nl,1}$ and $I_{nl,2}$ being the terms in $I_{nl}$ corresponding respectively in $\widetilde{\Theta}_T$ and in $B_r\times(0,T)$.

First, we compute $I_{nl,1}$ similarly to \cite{E,FI}. By Appendix, (\ref{psi}), (\ref{bounded}) and (\ref{alpha}), and $\widetilde{\Theta}=\Theta\cup\omega_0$, for $\lambda>2$ and $s>0$ large enough, one obtains
\begin{align}\label{I,nl1,4}
\nonumber I_{nl,1}
&\geq s^3\lambda^4\int_{\widetilde{\Theta}_T}(2-\alpha)p^2\theta^3\phi^3|\nabla\psi|^4z^2\,{\rm d}x\,{\rm d}t-Cs^3\lambda^3\int_{\widetilde{\Theta}_T}\theta^3\phi^3z^2\,{\rm d}x\,{\rm d}t\\ 
&\geq C_7s^3\lambda^4\int_{\widetilde{\Theta}_T}\theta^3\phi^3z^2\,{\rm d}x\,{\rm d}t-Cs^3\lambda^4\int_{\omega_{0,T}}\theta^3\phi^3z^2\,{\rm d}x\,{\rm d}t,
\end{align}
where,  $C_7=\frac{1}{2}p_1^2\delta^4$.

Next, we compute $I_{nl,2}$.
\begin{align}\label{I,nl2,1}
\nonumber I_{nl,2}&=s^3\int_{B_r\times(0,T)}(2-\alpha)p^2\theta^3(\lambda^2-2\lambda)\frac{|x|^{\lambda}}{r^{\lambda}}
\left(1+\lambda\frac{|x|^{\lambda-2}}{r^{\lambda}}\right)^2z^2\,{\rm d}x\,{\rm d}t\\ \nonumber
&+2s^3\int_{B_r\times(0,T)}p^2\theta^3\bigg(1+\lambda\frac{|x|^{\lambda-2}}{r^\lambda}\bigg)^3
|x|^2z^2\,{\rm d}x\,{\rm d}t\\ \nonumber
&+s^3\int_{B_r\times(0,T)}p\,\theta^3\left(1+\lambda\frac{|x|^{\lambda-2}}{r^\lambda}\right)^3
\left(\nabla p\cdot x\right)|x|^2z^2\,{\rm d}x\,{\rm d}t\\ \nonumber
&-s^3\int_{B_r\times(0,T)}n\,\alpha p^2\theta^3\bigg(1+\lambda\frac{|x|^{\lambda-2}}{r^\lambda}\bigg)^3
|x|^2z^2\,{\rm d}x\,{\rm d}t\\ 
&-s^2\int_{B_r\times(0,T)}\alpha^2p^2 \theta^2\left(n
+\lambda(\lambda+n-2)\frac{|x|^{\lambda-2}}{r^{\lambda}}\right)^2z^2\,{\rm d}x\,{\rm d}t.
\end{align}
Then for all $(x,t)\in B_r\times(0,T),$ and by the definition of $\alpha,$ taking $s_1(\lambda)>0$ large enough,, for instance $s_1(\lambda)> C\lambda^4$, we can obtain that
\begin{align}\label{2.59}
&\alpha^2 s^2p^2\theta^2\left(n
+\lambda(\lambda+n-2)\frac{|x|^{\lambda-2}}{r^{\lambda}}\right)^2\leq \frac{n}{3}\alpha\, s^3p^2\theta^3\left(1+\lambda\frac{|x|^{\lambda-2}}{r^\lambda}\right)^3
|x|^2,
\end{align}
for $\forall s>s_1(\lambda).$
Combining (\ref{a}), (\ref{r}), (\ref{alpha}),  (\ref{I,nl2,1}), (\ref{2.59}),
 there exist a constant $C_8>0$ such that for $\lambda$ and $s$ large enough, one has
\begin{align}\label{I,nl2,2}
\nonumber
I_{nl,2}&\geq s^3\int_{B_r\times(0, T)}\theta^3\left(1+\lambda\frac{|x|^{\lambda-2}}{r^{\lambda}}\right)^3p\left\{p\left(2-n\alpha -\frac{n\alpha}{3}\right)-|\nabla p\cdot x|\right\}|x|^2z^2{\rm d}x{\rm d}t\\ \nonumber
&\geq s^3\int_{B_r\times(0, T)}\theta^3\left(1+\lambda\frac{|x|^{\lambda-2}}{r^{\lambda}}\right)^3p\left(\frac{2}{3}p-p_3r\right)|x|^2z^2{\rm d}x{\rm d}t\\ \nonumber
&\geq s^3\int_{B_r\times(0, T)}\theta^3\left(1+\lambda\frac{|x|^{\lambda-2}}{r^{\lambda}}\right)^3p\left(\frac{2}{3}p_1-p_3r\right)|x|^2z^2{\rm d}x{\rm d}t\\
&\geq C_8s^3 \int_{B_r\times(0,T)}\theta^3|x|^2z^2\,{\rm d}x\,{\rm d}t.
\end{align}
Here, $C_8=\frac{2}{3}p_2\frac{\mu}{\mu^*}$.

Combining (\ref{I,nl1,4}) and (\ref{I,nl2,2}), we completed the proof of Lemma 3.4.

{\it Proof of Lemma 3.5.}

Let $\lambda>2$ be large enough. We note  that
\begin{gather}\label{t}
|\theta\theta'|\leq C\theta^3,\quad |\theta|^2\leq C\theta^3, \quad |\theta'|\leq C\theta^3,\quad |\theta''|\leq C\theta^{1+2/{k_0}},\qquad \forall t\in(0,T).
\end{gather}
By (\ref{a}), (\ref{t}), Appendix and the definition of $\sigma$, $\psi$ and $\xi$,   one can obtain
\begin{align}\label{dt1}
(x,t)\in\widetilde{\Theta}_T,\quad \vert{\partial_t\sigma}\vert\leq  s{\rm e}^{2\lambda\sup \psi}\vert\theta'\vert,\quad \vert{\partial_j\partial_t\sigma}\vert\leq Cs|\theta'|\lambda \phi
\end{align}
and
\begin{align}\label{dt2}
&(x,t)\in B_r\times(0,T),\quad \vert{\partial_t\sigma}\vert\leq  s{\rm e}^{2\lambda\sup \psi}\vert\theta'\vert,\quad \vert{\partial_j\partial_t\sigma}\vert\leq Cs|\theta'|\lambda|x|^2.
\end{align}

We can follow the proof of Lemma V.4 in \cite{J}. Similarly to Section V.3 in \cite{J}, by (1.2), we have
\begin{align*}
&|I_r|\leq C\left(s^2\lambda^2+s^2\lambda^2{\rm e}^{2\lambda\sup\psi}+se^{2\lambda\sup{\psi}}\frac{1}{\beta^q}\right)\int_{\Omega_T}\theta^3|x|^2z^2\,{\rm d}x\,{\rm d}t\\
&+C\left(s^2\lambda^2{\rm e}^{2\lambda \sup{\psi}}+s^2\lambda^2\right)\int_{\widetilde{\Theta}_T}\theta^3\phi^3 z^2\,{\rm d}x\,{\rm d}t
+C_{10}s{\rm e}^{2\lambda\sup{\psi}}\beta^{q'}\int_{\Omega_T}\theta\frac{z^2}{|x|^\gamma}\,{\rm d}x\,{\rm d}t
\end{align*}
for some constant $\beta>0$ , where 
$$
q=\frac{2+
\gamma}{\gamma}, \qquad q'=\frac{2+\gamma}{2}.
$$ For any given $\lambda>0$, we can choose $\beta=\beta(\lambda)$ such that
$$C_{10}{\rm e}^{2\lambda\sup{\psi}}\beta^{q'}=\frac{C_2}{2}.$$
Hence we deduce
$$\frac{1}{\beta^q}=C{\rm e}^{(4\lambda\sup{\psi})/\gamma} $$
for some constant $C>0$. It followings that
\begin{align*}
|I_r|&\leq C_{11}\left(s^2\lambda^2+s^2\lambda^2{\rm e}^{2\lambda\sup\psi}+s{\rm e}^{\lambda\sup{\psi}(2+4/\gamma)}\right)\int_{\Omega_T}\theta^3|x|^2z^2\,{\rm d}x\,{\rm d}t\\
&+C_{12}\left(s^2\lambda^2{\rm e}^{2\lambda \sup{\psi}}+s^2\lambda^2\right)\int_{\widetilde{\Theta}_T}\theta^3\phi^3 z^2\,{\rm d}x\,{\rm d}t
+\frac{C_2}{2}s\int_{\Omega_T}\theta\frac{z^2}{|x|^\gamma}\,{\rm d}x\,{\rm d}t.
\end{align*}
For any $\lambda>0$, there exists $s_0(\lambda)>s_1(\lambda)$ such that for all $s\geq s_0(\lambda)$, one has
$$C_{11}\left(s^2\lambda^2+s^2\lambda^2{\rm e}^{2\lambda\sup\psi}+s{\rm e}^{\lambda\sup{\psi}(2+4/\gamma)}\right)\leq \frac{C_8}{2}s^3 $$
and
$$C_{12}\left(s^2\lambda^2{\rm e}^{2\lambda \sup{\psi}}+s^2\lambda^2\right)\leq\frac{C_7}{4}s^3\lambda^4. $$

Therefore we obtain the inequality stated in Lemma 3.5.

{\it Proof of Lemma 3.7.}

Let us recall that $\omega$ and $\omega_0$ satisfy $\omega_0\subset\subset\omega$ and let us consider a smooth function $\xi:\overline{\Omega}\rightarrow\mathbb{R}$ such that
\begin{equation}\label{xi}
 \begin{cases}
 0\leq\xi(x)\leq1, \quad x\in\overline{\Omega},\\
 \xi(x)=1, \qquad x\in\omega_0,\\
 \xi(x)=0,\qquad x\in\overline{\Omega}\backslash\widetilde{\omega},
 \end{cases}
 \end{equation}
 where $\omega_0\subset\subset\widetilde{\omega}\subset\subset\omega$.
 By  (\ref{condition}), (\ref{theta}) and (\ref{pm}), one get
\begin{align*}
0&=\int_0^T\frac{\rm d}{\rm{d}t}\int_\omega\theta\phi\,\xi^2{\rm e}^{-2\sigma}w^2\,{\rm d}x\,{\rm d}t\\
&=-2\int_{\omega_T}\theta\phi\,\xi^2\partial_t\sigma{\rm e}^{-2\sigma}w^2\,{\rm d}x\,{\rm d}t
+2\int_{\omega_T}\theta\phi\,\xi^2{\rm e}^{-2\sigma}ww_t\,{\rm d}x\,{\rm d}t+\int_{\omega_T}\theta'\phi\,\xi^2{\rm e}^{-2\sigma}w^2\,{\rm d}x\,{\rm d}t\\
&=-2\int_{\omega_T}\theta\phi\,\xi^2\partial_t\sigma{\rm e}^{-2\sigma}w^2\,{\rm d}x\,{\rm d}t\mp\bigg(2\int_{\omega_T}\theta\phi\,\xi^2{\rm e}^{-2\sigma}w\,{\rm div}(p\nabla w)\,{\rm d}x\,{\rm d}t\\
&+2\int_{\omega_T}\theta\phi\,\xi^2{\rm e}^{-2\sigma}\mu\frac{w^2}{|x|^2}\,{\rm d}x\,{\rm d}t-2\int_{\omega_T}\theta\phi\,\xi^2{\rm e}^{-2\sigma}wg\,{\rm d}x\,{\rm d}t\bigg)
+\int_{\omega_T}\theta'\phi\,\xi^2{\rm e}^{-2\sigma}w^2\,{\rm d}x\,{\rm d}t.
\end{align*}
It follows that
\begin{align*}
&2\int_{\omega_T}\theta\phi\,p\,\xi^2 {\rm e}^{-2\sigma}|\nabla w|^2\,{\rm d}x\,{\rm d}t\\ \nonumber
&=\pm2\int_{\omega_T}\theta\phi\,\xi^2\partial_t\sigma{\rm e}^{-2\sigma}w^2\,{\rm d}x\,{\rm d}t\mp\int_{\Omega_T}\theta'\phi\,\xi^2{\rm e}^{-2\sigma}w^2\,{\rm d}x\,{\rm d}t
-4\int_{\omega_T}\theta\phi\,p\,\xi {\rm e}^{-2\sigma}w(\nabla\xi\cdot\nabla w)\,{\rm d}x\,{\rm d}t\\ \nonumber
&+4\int_{\omega_T}\theta\phi\,p\,\xi^2 {\rm e}^{-2\sigma}w(\nabla\sigma\cdot\nabla w)\,{\rm d}x\,{\rm d}t
+2\int_{\omega_T}\theta\phi\,\xi^2{\rm e}^{-2\sigma}\mu\frac{w^2}{|x|^2}\,{\rm d}x\,{\rm d}t
-2\int_{\omega_T}\theta\phi\,\xi^2{\rm e}^{-2\sigma}wg\,{\rm d}x\,{\rm d}t\\ \nonumber
&-2\lambda\int_{\omega_T}\theta\phi\,p\,\xi^2 {\rm e}^{-2\sigma}w(\nabla\psi\cdot\nabla w)\,{\rm d}x\,{\rm d}t.
\end{align*}
By (\ref{Ome}), (\ref{a}),  (\ref{Br}) and (\ref{2.27}),   there exist some positive constants $C_{13}, C_{14}, C_{15}$ (independent of $s$ and $\lambda$) and $\epsilon_1, \epsilon_2, \epsilon_3$ such that the following inequalities hold:
\begin{align*}
&\left\vert{-4\theta\phi\,p\,\xi {\rm e}^{-2\sigma}w(\nabla\xi\cdot\nabla w)}\right\vert\leq 2C_{13}\theta\phi\,\xi {\rm e}^{-2\sigma}|w|\,|\nabla w|
\leq C_{13}\left(\frac{1}{\epsilon_1}\theta\phi {\rm e}^{-2\sigma}w^2
+\epsilon_1\,\theta\phi\xi^2{\rm e}^{-2\sigma}|\nabla w|^2\right),\\
&\left\vert{4\theta\phi\,p\,\xi^2 {\rm e}^{-2\sigma}w(\nabla\sigma\cdot\nabla w)}\right\vert\leq 2C_{14}s\lambda\theta^2\phi^2\xi^2{\rm e}^{-2\sigma}|w|\,|\nabla w|\\
&\,\,\,\,\,\,\,\,\,\,\,\,\,\,\,\,\,\,\,\,\,\,\,\,\,\,\,\,\,\,\,\,\,\,\,\,\,\,\,\,\,\,\,\,\,\,\,\,\,\,\,\,\,\,\,\,\,\,\,\,\,\,\,\,\,\,\,\,\, \leq C_{14}\left(\epsilon_2\,\theta\phi\xi^2{\rm e}^{-2\sigma}|\nabla w|^2
+\frac{1}{\epsilon_2}s^2\lambda^2\phi^3\theta^3\xi^2{\rm e}^{-2\sigma}w^2\right),\\
&\left\vert{-2\theta\phi\,\xi^2{\rm e}^{-2\sigma}wg}\right\vert\leq \frac{1}{s\lambda^2}\xi^2{\rm e}^{-2\sigma}g^2+s\lambda^2\xi^2\theta^2\phi^2{\rm e}^{-2\sigma}w^2,\\
&\left\vert{-2\lambda\theta\phi\,p\,\xi^2 {\rm e}^{-2\sigma}w(\nabla\psi\cdot\nabla w)}\right\vert\leq 2C_{15}\lambda\theta\phi\,\xi^2 {\rm e}^{-2\sigma}|w|\,|\nabla w|\\
&\,\,\,\,\,\,\,\,\,\,\,\,\,\,\,\,\,\,\,\,\,\,\,\,\,\,\,\,\,\,\,\,\,\,\,\,\,\,\,\,\,\,\,\,\,\,\,\,\,\,\,\,\,\,\,\,\,\,\,\,\,\,\,\,\,\,\,\,\,\,\,\,\,\,\,\,  \leq C_{15}\left(\epsilon_3\,\theta\phi\xi^2{\rm e}^{-2\sigma}|\nabla w|^2+\frac{1}{\epsilon_3}\lambda^2\theta\phi \xi^2{\rm e}^{-2\sigma}w^2\right).
\end{align*}
Let $\epsilon_1, \epsilon_2, \epsilon_3$ satisfy : $p_1=C_{13}\epsilon_1+C_{14}\epsilon_2+C_{15}\epsilon_3$, where $p_1$ is as in (\ref{a}).
Therefore, by (\ref{t}) and (\ref{dt1}), for   $\forall \mu\leq\,p_1\mu^*$, one can obtain
\begin{align*}
\int_{\omega_0\times(0,T)}&\theta \phi {\rm e}^{-2\sigma}|\nabla w|^2\,{\rm d}x\,{\rm d}t\leq
\int_{\omega\times(0,T)}\theta\phi\,\xi^2 {\rm e}^{-2\sigma}|\nabla w|^2\,{\rm d}x\,{\rm d}t\\ \nonumber
&\leq C\int_{\omega\times(0,T)}\theta\phi \xi^2|\partial_t\sigma|{\rm e}^{-2\sigma}w^2\,{\rm d}x\,{\rm d}t+C\int_{\omega\times(0,T)}|\theta'|\phi \xi^2{\rm e}^{-2\sigma}w^2\,{\rm d}x\,{\rm d}t\\
&+C\int_{\omega\times(0,T)}\theta\phi {\rm e}^{-2\sigma}w^2\,{\rm d}x\,{\rm d}t+Cs^2\lambda^2\int_{\omega\times(0,T)}\theta^3\phi^3 \xi^2{\rm e}^{-2\sigma}w^2\,{\rm d}x\,{\rm d}t\\
&+2\int_{\omega\times(0,T)}\theta\phi \xi^2{\rm e}^{-2\sigma}\frac{\mu}{|x|^2}w^2\,{\rm d}x\,{\rm d}t
+\frac{1}{s\lambda^2}\int_{\omega\times(0,T)}\xi^2{\rm e}^{-2\sigma}g^2\,{\rm d}x\,{\rm d}t\\
&+s\lambda^2\int_{\omega\times(0,T)}\theta^2\phi^2\xi^2{\rm e}^{-2\sigma}w^2\,{\rm d}x\,{\rm d}t+C\lambda^2\int_{\omega\times(0,T)}\theta\phi\xi^2{\rm e}^{-2\sigma}w^2\,{\rm d}x\,{\rm d}t\\
&\leq C(s{\rm e}^{2\lambda\sup\psi}+1+s^2\lambda^2+s\lambda^2+\lambda^2)\int_{\omega\times(0,T)}\theta^3\phi^3{\rm e}^{-2\sigma}w^2\,{\rm d}x\,{\rm d}t\\
&+\frac{1}{s\lambda^2}\int_{\omega\times(0,T)}{\rm e}^{-2\sigma}g^2\,{\rm d}x\,{\rm d}t.
\end{align*}
Then, for $\lambda$ large enough,  we get (\ref{2.26}).

Note that all of the constants in the proof of Lemma 3.7 are independent of $\mu,s,\lambda$.

This ends the proof of Lemma 3.7.

 $\hfill\square$ $\vspace{2mm}$

\section{Non Uniform Stabilization in the Case $\mu>p_2\mu^*$ }
\setcounter{equation}{0}
As in \cite{E}, following the ideas of optimal control, for any $u_0\in L^2(\Omega),$ we consider the functional
\begin{gather}\label{1.10}
J_{u_0}(u,f)=\frac{1}{2}\int_{\Omega_T}|u(t,x)|^2\,{\rm d}x\,{\rm d}t+\frac{1}{2}\int_0^T\Vert{f(t)}\Vert_{H^{-1}(\Omega)}^2\,{\rm d}t,
\end{gather}
defined on the set
\begin{align}
\mathbf{U}(u_0)=\{&(u,f)\in L^2((0,T); H_0^1(\Omega))\times L^2((0,T);H^{-1}(\Omega))\,\rm{ such\,\,  that}\\ \nonumber
 &u \,\rm {satisfies \,\, (1.1) \,\,with \,\, f \,\,as\,\,  in\,\, (\ref{*})}\}.
\end{align}
We say that we can stabilize system (1.1) if we can find a constant $C$ such that
\begin{gather}\label{1.12}
\forall u_0\in L^2(\Omega),\qquad \inf\limits_{(u,f)\in \mathbf{U}(u_0)} J_{u_0}(u,f)\leq C\Vert{u_0}\Vert^2_{L^2(\Omega)}.
\end{gather}
Of course, this property strongly depends on the set $\omega$ where the stabilization is effective. Especially, when $0\in \omega$, (\ref{1.12}) holds (see Section 5). When $0\notin \omega$, the situation is more intricate. Therefore we focus our study on this particular case, and give a severe obstruction, in this case, to the stabilization property (\ref{1.12}).

More precisely, for $\epsilon>0,$ we approximate (1.1) by the systems
\begin{equation}\label{1.13}
 \begin{cases}
 \partial_tu(x,t)-{\rm div}(p(x)\nabla u(x,t))-\frac{\mu}{|x|^2+\epsilon^2}u(x,t)=f(x,t), \quad{(x,t)\in\Omega_T,}\\
 u(x,t)=0, \qquad{ (x,t)\in\Sigma_T,}\\
 u(x,0)=u_0(x),\qquad {x\in\Omega.}
 \end{cases}
 \end{equation}
 For these approximate problems, the Cauchy problem is well-posed. Therefore we can consider the functionals
\begin{gather}
J_{u_0}^\epsilon(f)=\frac{1}{2}\int_{\Omega_T}|u(t,x)|^2\,{\rm d}x\,{\rm d}t+\frac{1}{2}\int_0^T\Vert{f(t)}\Vert_{H^{-1}(\Omega)}^2\,{\rm d}t,
\end{gather}
where $f\in L^2((0,T);H^{-1}(\Omega))$ is localized in $\omega$ in the sense of (\ref{*}) and $u$ is the corresponding solution of (\ref{1.13}). We prove the following:
\begin{theorem}
Assume that $\mu>p_2\mu^*$, and that $0\notin\bar{\omega}$. There is no constant $C$ such that for all $\epsilon>0$, and for all $u_0\in L^2(\Omega).$
\begin{equation}
\inf\limits_{\small{\begin{aligned}
&f \in L^2\left((0; T); H^{-1}(\Omega)\right),\\
&\,\,\,\,\,\,\,\,\,\,\,\,\,\,\,f\,\, as\,\,in\,\, (\ref{*})
\end{aligned}}}
J_{u_0}^\epsilon(f)\leq C\Vert{u_0}\Vert_{L^2(\Omega)}^2.
\end{equation}
\end{theorem}
In particular, this result implies that the stabilization of (1.1) is impossible to attain through regularization processes when $\mu> p_2\mu^*$ and $0\notin \bar{\omega}$, and that we cannot prevent the system from blowing up.

The goal of this section is to prove Theorem 4.1. The proof  will rely on the proof of Theorem 1.2  in \cite{E}. It is divided into two main steps.

First, we prove some basic estimates on the spectrum of the operator
\begin{gather}
L^\epsilon=-{\rm div}(p(x)\nabla\,)-\frac{\mu}{|x|^2+\epsilon^2}
\end{gather}
on $\Omega$ with Dirichlet boundary conditions, especially on the first eigenvalue $\lambda_0^\epsilon$ and the corresponding eigenfunction $\phi_0^\epsilon$. This will be done in Subsection 4.1.

Second, we verify Theorem 4.1 in Subsection 4.2 by giving a lower bound on the quantity $J_{\phi_0^\epsilon}^\epsilon$ that goes to infinity when $\epsilon\rightarrow0$.

\subsection{Spectral Estimates}
Since for $\epsilon>0$, the function $1/(|x|^2+\epsilon^2)$ is smooth and bounded in $\Omega$, the spectrum of $L^\varepsilon$ is formed by a sequence of real eigenvalues $\lambda_0^\epsilon\leq\lambda_1^\epsilon\leq\cdots\leq\lambda_n^\epsilon\leq\cdots$, with $\lambda_n^\epsilon\rightarrow+\infty$ when $n\rightarrow\infty$. The corresponding eigenvectors $\phi_n^\epsilon$ are a basis of $L^2(\Omega),$ orthonormal with respect to the $L^2$ scale product. We choose $\phi_n^\epsilon$ of unit $L^2-$ norm.

In the sequel, we focus on the bottom of the spectrum - the most explosive mode.
\begin{lemma}
Assume that $\mu>p_2\mu^*$. Then we have that
\begin{gather}\label{3.2}
\lim\limits_{\epsilon\rightarrow0}\lambda_0^\epsilon=-\infty,
\end{gather}
and for all $\tau>0,$
\begin{gather}\label{3.3}
\lim\limits_{\epsilon\rightarrow0}\Vert{\phi_0^\epsilon}\Vert_{H^1(\Omega\backslash \overline{{B}}(0,\tau))}=0.
\end{gather}
\end{lemma}
{\it Proof.}
We argue by contradiction, and assume that $\lambda_0^\epsilon$ is bounded from below for a subsequence by a real number $C$. Then, from the Rayleigh formula (e.g.,  Remark (ii) of Theorem 2 in Chapter 6.5 of \cite{LC}), we get
\begin{gather*}
\forall \epsilon>0, \quad u\in H_0^1(\Omega), \quad \mu\int_\Omega\frac{u^2}{|x|^2+\epsilon^2}\,{\rm d}x\leq
\int_{\Omega}p\,|\nabla u|^2\,{\rm d}x-C\int_\Omega u^2\,{\rm d}x.
\end{gather*}
Taking $u\in\wp(\Omega)$, we pass to the limit $\epsilon\rightarrow0$ and get
\begin{gather}\label{3.4}
\mu\int_\Omega\frac{u^2}{|x|^2}\,{\rm d}x\leq \int_{\Omega}p\,|\nabla u|^2\,{\rm d}x-C\int_\Omega u^2\,{\rm d}x,
\end{gather}
that must therefore hold for any $u\in H_0^1(\Omega)$ by a density argument.

Now, there exists $\tau_0>0$ such that $B(0,\tau_0)\subset\Omega$. We then choose $u\in H_0^1(B(0,\tau_0))$ that we extend by 0 on $\mathbb{R}^n$, and define for $a\geq1$
$$u_a(r)=a^nu(ar).$$
These functions are in $H_0^1(B(0,\tau_0))$, and therefore in $H_0^1(\Omega)$, and we apply (\ref{3.4}) to them:
$$a^2\left(\mu\int_{\Omega}\frac{u^2}{|x|^2}\,{\rm d}x-\int_{\Omega}p\,|\nabla u|^2\,{\rm d}x\right)\,{\rm d}x\leq-C\int_\Omega u^2\,{\rm d}x.$$
Passing to the limit $a\rightarrow\infty$, we obtain that
$$\forall u\in H_0^1(B(0,\tau_0)), \quad \mu\int_{\Omega}\frac{u^2}{|x|^2}\leq \int_{\Omega}p\,|\nabla u|^2\,{\rm d}x.$$
Therefore we should have that $\mu\leq p_2\mu^*$, since this is the Hardy inequality (\ref{1.3}) in the set $B(0,\tau_0)$, and then we have a contradiction.

Now, consider the first eigenvector $\phi_0^\epsilon\in H_0^1(\Omega)$ of $L^\epsilon$:
\begin{gather}\label{3.5}
-{\rm div}(p\,\nabla\phi_0^\epsilon)-\frac{\mu}{|x|^2+\epsilon^2}\phi_0^\epsilon=\lambda_0^\epsilon\phi_0^\epsilon, \quad \forall x\in \Omega.
\end{gather}
Remark that since the potential is smooth in $\Omega$, the function $\phi_0^\epsilon$ is smooth by classical elliptic estimates.

Set $\tau>0$. Let $\epsilon >0$ be sufficiently small such that $\lambda_0^\epsilon<0$. Let $\eta_\tau$ be a nonnegative and sufficiently smooth function that vanishes in $B(0,\tau/2)$ and equals 1 in $\mathbb{R}^n\backslash B(0,\tau)$ with $\Vert{\eta_\tau}\Vert_\infty\leq1$. Multiplying (\ref{3.5}) by $\eta_\tau\phi_0^\epsilon$ and integrating by parts, we get:
\begin{align}\label{3.6}
&\int_{\Omega}\eta_\tau \,p\,|\nabla\phi_0^\epsilon|^2\,{\rm d}x
+|\lambda_0^\epsilon|\int_{\Omega}\eta_\tau|\phi_0^\epsilon|^2\,{\rm d}x
=\mu\int_{\Omega}\eta_\tau\frac{|\phi_0^\epsilon|^2}{|x|^2+\epsilon^2}\,{\rm d}x\\ \nonumber
&+\frac{1}{2}\int_{\Omega}|\phi_0^\epsilon|^2p\,\triangle\eta_\tau\,{\rm d}x
+\frac{1}{2}\int_{\Omega}|\phi_0^\epsilon|^2(\nabla\eta_\tau\cdot\nabla p)\,{\rm d}x.
\end{align}
Therefore, since $\phi_0^\epsilon$ is of unit $L^2-$ norm, due to the particular choice of $\eta_\tau$, we get
\begin{align*}
|\lambda_0^\epsilon|\int_{\Omega\backslash B(0,\tau)}|\phi_0^\epsilon|^2\,{\rm d}x
&\leq \frac{4\mu}{\tau^2}+C\Vert{\eta_\tau}\Vert_{C^2(\Omega)}.
\end{align*}
Since $|\lambda_0^\epsilon|\rightarrow\infty$ when $\epsilon\rightarrow0$, we get that for any $\tau>0$,
\begin{gather}\label{3.7}
\lim_{\epsilon\rightarrow0}\int_{\Omega\backslash B(0,\tau)}|\phi_0^\epsilon|^2\,{\rm d}x=0.
\end{gather}
Besides, still using (\ref{3.6}) and the particular form of $\eta_\tau$
\begin{align*}
\int_{\Omega\backslash B(0,\tau)} |\nabla\phi_0^\epsilon|^2\,{\rm d}x\leq
 C\left(\frac{4\mu}{\tau^2}+\Vert{\eta_\tau}\Vert_{C^2(\Omega)}\right)\int_{\Omega\backslash B(0,\tau/2)}|\phi_0^\epsilon|^2\,{\rm d}x.
\end{align*}
Therefore the proof of (\ref{3.3}) is completed by using (\ref{3.7}) for $\tau/2$ instead of $\tau$.

\subsection{Proof of Theorem 4.1}
Fix $\epsilon>0$, and choose $u_0^\epsilon$, which is of unit $L^2$-norm. Our goal is to prove that
\begin{equation}
\inf\limits_{\small{\begin{aligned}
& {f \in L^2\left((0; T); H^{-1}(\Omega)\right),}\\
 &\,\,\,\,\,\,\,\,\,\,\,\,\,\,f\,\, as\,\,in\,\, (\ref{*})
 \end{aligned}}}
 J_{u_0}^\epsilon(f) \rightarrow \infty ,\qquad \epsilon \to 0.  
\end{equation}
Let $f \in L^2\left((0; T); H^{-1}(\Omega)\right),
f$ as in (\ref{*}), and consider $u$ the corresponding solution of (\ref{1.13}) with initial data $u_0=\phi_0^\epsilon$.

Set 
$$a(t)=\int_\Omega u(t,x)\phi_0^\epsilon(x)\,{\rm d}x,\qquad b(t)=\left<f(t),\phi_0^\epsilon\right>_{H^{-1}(\Omega)\times H_0^1(\Omega)}.$$
Then $a(t)$ satisfies the  equation
$$a'(t)+\lambda_0^\epsilon a(t)=b(t), \qquad a(0)=1.$$
We refer to Section 3.2 in \cite{E}, so that 
 for any $f \in L^2\left((0; T); H^{-1}(\Omega)\right),
f$ as in (\ref{*}), we get
$$J_{u_0}^\epsilon(f)\geq \inf\left\{  \frac{{\rm e}^{2|\lambda_0^\epsilon|T}-1}{16|\lambda_0^\epsilon|},
\frac{|\lambda_0^\epsilon|}{4\Vert{\phi_0^\epsilon}\Vert^2_{H^1(\omega)}} \left(1-{\rm e}^{-2|\lambda_0^\epsilon|T} \right) \right\}.$$ 
 This bound blows up when $\epsilon\to 0$ from the estimates (\ref{3.2}) and (\ref{3.3}). Indeed, since $0\notin\bar\omega,$ we can choose $\tau>0$ small enough such that $\omega\subset\Omega\backslash B(0,\tau)$ and therefore
 $$ \Vert{\phi_0^\epsilon}\Vert^2_{H^1(\omega)} \leq \Vert{\phi_0^\epsilon}\Vert_{H^1(\Omega\backslash B(0,\tau))}\rightarrow 0, \quad \epsilon\to 0.  $$

\section{Conclusions}
In this article we studied a parabolic equation with variable coefficients in the principal part, which have an inverse-square potential $-\mu/|x|^2$, from a control point of view, in the case $0\leq\mu < \frac{p_1^2}{p_2}\mu^*$, and the case $\mu>p_2\mu^*$.

When $0\leq\mu< \frac{p_1^2}{p_2}\mu^*$, we have addressed the null-controllability problem for a distributed control in an arbitrary open subset of $\Omega$. To this end, we have derived a  Carleman inequality (\ref{CE}) inspired by the articles \cite{E,FI,J} and \cite{VZ1}. When $\mu>p_2\mu^*$, we have shown that we cannot uniformly stabilize regularized approximations of system (1.1) with a control supported in $\omega$ when $0\notin\omega$ as in \cite{E}.

To complete this result, as in \cite{E}, we comment the case $0\in\omega$, for which the stabilization property (\ref{1.10}) holds. Given $u_0\in L^2(\Omega)$, we claim that we can find $u\in L^2((0,T); H_0^1(\Omega))$ and $f\in L^2((0,T); H^{-1}(\Omega))$ as in (\ref{*}) such that $u$ is the solution of (1.1) and that $J_{u_0}(u,f)\leq C\Vert{u_0}\Vert^2_{L^2(\Omega)}$.

Indeed, as in \cite{E}, by $\chi$ a smooth function that equals 1 in a neighborhood of 0 and vanishing outside $\omega$. Then consider the solution $u$ of
\begin{equation*}
 \begin{cases}
 \partial_tu(x,t)-{\rm div}(p\nabla u(x,t))-(1-\chi)\frac{\mu}{|x|^2}u(x,t)=0, \quad{(x,t)\in\Omega_T,}\\
 u(x,t)=0,\qquad{ (x,t)\in\Sigma_T,}\\
 u(x,0)=u_0(x),\qquad {x\in\Omega,}
 \end{cases}
 \end{equation*}
 which satisfies $u\in L^2((0,T); H_0^1(\Omega)),$ and $\Vert{u}\Vert_{L^2(0,T; H_0^1(\Omega))}\leq C\Vert{u_0}\Vert^2_{L^2}$ for some constant $C$. Then taking $f=\chi\mu u/{|x|^2}\in L^2((0,T); H^{-1}(\Omega))$ provides an admissible stabilizer with the required property (\ref{*}).

 The same argument can also be applied to derive the null controllability property for (1.1) when $0\in\omega$. Indeed, the  results in \cite{FI} proves that there exists a control $v\in L^2((0,T)\times \omega)$ such that the solution of
\begin{equation*}
 \begin{cases}
 \partial_tu(x,t)-{\rm div}(p\nabla u(x,t))-(1-\chi)\frac{\mu}{|x|^2}u(x,t)=v, \quad{(x,t)\in\Omega_T,}\\
 u(x,t)=0,\qquad{ (x,t)\in\Sigma_T,}\\
 u(x,0)=u_0(x),\qquad {x\in\Omega.}
 \end{cases}
 \end{equation*}
 satisfies $u(T)=0$. Besides, the norms of $v$ in $L^2((0,T)\times\omega)$ and $u$ in $L^2((0,T); H_0^1(\Omega))$ are bounded by the norm of $u_0$
in $L^2(\Omega)$. Then taking $f=v+\mu \chi u/{|x|^2}$ provides a control in $L^2((0,T); H^{-1}(\Omega))$ for (1.1) that drives the solution to 0 in time $T$.

\vspace{5mm}
{\bf Acknowledgements}

 The first author is supported by the  scholarship from China Scholarship Council (CSC).

\section*{Appendix}
By (\ref{sigma}) and direct calculations, we have
\begin{gather*}
\partial_t\sigma=s\theta'({\rm e}^{2\lambda\sup \psi}-\frac{|x|^2}{2}-{\rm e}^{\lambda\psi}),\quad \partial_k\sigma=-s\theta(x_k+\lambda\phi\partial_k\psi),\\ \nonumber
\partial_j\partial_k\sigma=-s\theta(\delta_{jk}+\lambda\phi\partial_j\partial_k\psi+\lambda^2\phi(\partial_j\psi)\partial_k\psi),\\ \nonumber
\partial_j\partial_k\partial_m\sigma=-s\theta\{\lambda\phi\partial_j\partial_k\partial_m\psi+\lambda^2\phi(\partial_j\psi)\partial_k\partial_m\psi\\ \nonumber
+\lambda^2\phi(\partial_j\partial_k\psi)\partial_m\psi+\lambda^2\phi(\partial_j\partial_m\psi)\partial_k\psi
+\lambda^3\phi(\partial_j\psi)(\partial_k\psi)\partial_m\psi\},\\ \nonumber
\sum_{j,k=1}^n\partial^2_j\partial_k^2\sigma=-s\theta\bigg\{\sum_{j,k=1}^n\
\bigg(\lambda\phi\partial^2_j\partial_k^2\psi
+2\lambda^2\phi(\partial_j\partial_k\psi)^2+4\lambda^2\phi(\partial^2_j\partial_k\psi)\partial_k\psi
+4\lambda^3\phi(\partial_j\psi)(\partial_k\psi)\partial_k\partial_j\psi\bigg)\\
+\lambda^4\phi|\nabla\psi|^4+\lambda^2\phi(\triangle\psi)^2+2\lambda^3\phi|\nabla \psi|^2\triangle\psi\bigg\}.\\ \nonumber
\end{gather*}
In $B_r$, we have
\begin{gather*}
\quad \phi=\frac{|x|^\lambda}{r^\lambda},\quad \partial_k\psi=\frac{x_k}{|x|^2},\quad
\partial_j\partial_k\psi=\frac{\delta_{kj}}{|x|^2}-2\frac{x_kx_j}{|x|^4},\\ \nonumber
\partial_k\partial_j\partial_m\psi=-2\frac{\delta_{mj}x_k}{|x|^4}-2\frac{\delta_{km}x_j}{|x|^4}-2\frac{\delta_{kj}x_m}{|x|^4}
+8\frac{x_mx_jx_k}{|x|^6},\\ \nonumber
\sum_{j,k=1}^n\partial_j^2\partial_k^2\psi=(-2n^2+12n-16)\frac{1}{|x|^4}.
\end{gather*}

\renewcommand{\baselinestretch}{0.3}

\end{document}